\documentclass[11pt,a4paper]{amsart}
\usepackage{amssymb}
\usepackage{epsfig}

\begin{document}
\newtheorem{thrm}{Theorem}
\newtheorem{thmf}{Th{\'e}or{\`e}me}
\newtheorem{thmm}{Theorem}
\def\thethmm{\ref{twist}$'$}
\newenvironment{thint}[1]{{\flushleft\sc{Th{\'e}or{\`e}me}}
      {#1}. \it}{\medskip} 
\newenvironment{thrm.}{{\flushleft\bf{Theorem}}. \it}{\medskip} 
\newenvironment{thrme}[1]{{\flushleft\sc{Theorem}}
      {#1}. \it}{\medskip} 
\newenvironment{propint}[1]{{\flushleft\sc{Proposition}}
      {#1}. \it}{\medskip} 
\newenvironment{corint}[1]{{\flushleft\sc{Corrolaire}}
      {#1}. \it}{\medskip} 
\newtheorem{corr}{Corollary}
\newtheorem{corf}{Corollaire}
\newtheorem{prop}{Proposition}
\newtheorem{defi}{Definition}
\newtheorem{deff}{D{\'e}finition}
\newtheorem{lem}{Lemma}
\newtheorem{lemf}{Lemme}
\def\fy{\varphi}
\def\ul{\underline}
\def\obsf{{\flushleft\bf Remarque. }}
\def\obs{{\flushleft\bf Remark. }}
\def\R{\mathbb{R}}
\def\C{\mathbb{C}}
\def\Z{\mathbb{Z}}
\def\N{\mathbb{N}}
\def\P{\mathbb{P}}
\def\Q{\mathbb{Q}}
\def\p{\pi_1(X)}
\def\Re{\mathrm{Re}}
\def\Im{\mathrm{Im}}
\def\H{\mathbf{H}}
\def\Hs{{Nil}^3}
\def\L{L_\alpha}
\def\h{\frac{1}{2}}
\def\M{\P(E)\times\P(E)^*\smallsetminus\mathcal{F}}
\def\m{\cp2\times{\cp2}^*\smallsetminus\mathcal{F}}
\def\g{\mathfrak{g}}
\def\l{\mathcal{L}}
\def\k{\mathfrak{h}}
\def\s3{\mathfrak{s}^3}
\def\nb{\nabla}
\def\Re{\mathrm{Re}}
\def\Im{\mathrm{Im}}
\newcommand\cp[1]{\mathbb{CP}^{#1}}
\renewcommand\o[1]{\mathcal{O}({#1})}
\renewcommand\d[1]{\partial_{#1}}
\def\sq{\square_X}
\def\gp{\dot\gamma}
\def\j{\mathcal{J}}
\def\iif{\mbox{\bf\em{I}}f}

\title{On the Weyl tensor of a self-dual complex 4-manifold}
\author{Florin Alexandru Belgun}
\thanks{AMS classification~: 53C21, 53A30, 32C10.}
\date{February 9, 1999.}

\begin{abstract}
We study complex 4-manifolds with holomorphic self-dual conformal structures,
and we obtain an interpretation of the Weyl tensor of such a manifold
as the projective curvature of a field of cones on the ambitwistor
space. In particular, its vanishing is implied by the existence of
some compact, simply-connected, null-geodesics. We also relate the
Cotton-York tensor of an umbilic hypersurface to the Weyl tensor of
the ambient. As a consequence, a conformal 3-manifold or a self-dual
4-manifold admitting a rational curve as a null-geodesic is conformally flat. We
show that the projective structure of the $\beta$-surfaces of a
self-dual manifold is flat.
\end{abstract}

\maketitle

\section{Introduction}\label{s.w1}
Twistor theory, created by Penrose \cite{p1}, establishes a close
relationship between conformal Riemannian geometry in dimension 4, and
(almost) complex geometry in dimension 3. In particular, to a
Riemannian manifold $M$ for which the part $W^-$ of the Weyl tensor
vanishes identically ({\it self-dual}), one associates its {\it
  twistor space} $Z$, a complex 3-manifold containing rational curves
with normal bundle $\o1\oplus\o1$, and admitting a real structure with
no fixed points \cite{ahs}, \cite{hi}, \cite{besse}.

The space of such curves is a complex 4-manifold $M^\C$ \cite{kod}
with a holomorphic conformal structure and is, therefore, a conformal
complexification of $M$ \cite{ahs}, \cite{hi}, \cite{besse}.

As the conformal geometry of $M$ is encoded by the complex geometry of
$Z$, we ask ourselves what holomorphic object on $Z$ corresponds to
$W^+$, the Weyl tensor of the self-dual manifold $M$. It seems that
this question, although natural, has not been considered in the
literature, and maybe a reason for that is that the answer appears to
be a highly non-linear object.

This object is more easily understood in the framework of
complex-Rie\-mannian geometry (see Section \ref{s.w2}): following LeBrun
\cite{leb4}, we (locally) introduce the space $B$ of complex null-geodesics of {\bf M}
({\it ambitwistor space}). 
For a self-dual (complex) 4-manifold {\bf M},
its (local) twistor space is then defined as the 3-manifold of
$\beta$-surfaces (some totally geodesic isotropic surfaces, see
Section \ref{s.w2}).

The ambitwistor space $B$ and (in the self-dual case) the twistor
space $Z$ completely describe the conformal structure of {\bf M}. In
particular, a null-geodesic $\gamma$ in {\bf M} corresponds to the set of rational
curves in $Z$ tangent to a 2-plane. The union of these curves, called
the {\it integral $\alpha$-cone} of $\gamma$ (see Section \ref{s.w3}), is
lifted to a (linearized) $\alpha$-cone in $T_\gamma B$. Our first
result (Theorem \ref{t1}) is that the Weyl tensor of {\bf M} is equivalent
to the projective curvature (see Section \ref{s.w4}) of the field of
$\alpha$-cones on $B$. In particular, if such a cone is {\it flat},
then $W^+$ vanishes on certain isotropic planes in {\bf M}.
\medskip

We use Theorem \ref{t1} to investigate global properties of a
self-dual manifold {\bf M}: If the integral $\alpha$-cone of $\gamma$ is
part of a smooth surface in $Z$, then the linearized $\alpha$-cone is
flat (Theorems \ref{twist}, \ref{3'}). In particular, the space
$\mathbf{M}_0$ of rational 
curves of $Z$ with normal bundle $\o1\oplus\o1$ is compact iff
$Z\simeq\cp3$. On the other hand, it is known, from a theorem of
Campana \cite{camp}, that, for a compact twistor space $Z$, $\mathbf{M}_0$ can
be compactified within the space of analytic cycles iff $Z$ is
Moishezon. It appears then that the conformal structure does not
extend smoothly to the compactification.

A good illustration of what happens in the non-flat (self-dual) case
is the K{\"a}hler-Einstein manifold $\cp2$ whose twistor space is known
to be the manifold of flags in $\C^3$ \cite{ahs}, see Section \ref{s.w8}.

Different methods allow us to generalize Theorem 3 to non
geodesically-connec\-ted self-dual manifolds: We show (Theorem
\ref{compact}) that if a self-dual manifold admits a compact,
simply-connected, null-geodesic, then it is conformally flat. 
We also note that
the rational curves in $Z$, corresponding to the points of {\bf M} (see
Section \ref{s.w2}) are then geodesics of some {\it projective structure} of $Z$ iff
the latter is projectively flat (Corollary \ref{twfl}).

The isotropic, totally geodesic surfaces (called $\beta$-surfaces) in
a self-dual manifold {\bf M} 
appear to have a projective structure, given by the null-geodesics of {\bf M}
contained in it (Section \ref{s.w6}). We show that it is {\it flat}
(i.e. locally equivalent to $\cp2$) (Corollary \ref{prflat}), and we
obtain a classification of the compact $\beta$-surfaces of a self-dual
4-manifold (Theorem \ref{clasb}).
\smallskip

Theorem \ref{compact} can be adapted for conformal 3-manifolds : A
conformal 3-manifold admitting a rational curve as a 
null-geodesic is conformally flat (Theorem
\ref{comp3}). The conformal geometry in dimensions 3 and 4 are
related, as any geodesically convex 3-manifold $Q$
can be realized as the {\it conformal infinity} of a self-dual 
4-manifold {\bf M} \cite{leb1}. 
In particular, $Q$ is umbilic in {\bf M}, and
we relate, in Section \ref{s.w7}, the conformal invariants of the 2
manifolds: the Cotton-York tensor of $Q$ is identified to the
derivative, in the normal direction, of the Weyl tensor of {\bf M}
(Theorem \ref{umbilic}).  This result can be equally stated in the
real framework.
\medskip

The paper is organized as follows: in Section \ref{s.w2} we recall the
classical results of the twistor theory (especially for complex
4-manifolds), in Section \ref{s.w3} we introduce the $\alpha$-cones on the
(ambi-)twistor space, and, in Section \ref{s.w4}, we prove the equivalence
between the projective curvature of the latter and the Weyl tensor
$W^+$ of {\bf M}. Section \ref{s.w5} is devoted to the proof of some results of the
type ``compactness implies conformal (projective) flatness'': Theorems
\ref{twist}, \ref{3'} and \ref{compact}, mentioned above. We study the projective
structure of $\beta$-surfaces in Section \ref{s.w7}, and we illustrate the above results on
the special case of the self-dual manifold $\cp2$, in Section \ref{s.w8}.

\bigskip

{\noindent\sc Acknowledgements } The author is deeply indebted to Paul
Gauduchon, for his care in reading the manuscript and for his constant
help during the research and redaction.

\section{Preliminaries}\label{s.w2}

The content of this paper makes use of {\it complex-Riemannian
  geometry} (with the exception of Theorem \ref{umbilic} and Corollary
  \ref{cumb}, which
  hold also in the real framework). Complex-Riemannian geometry is
  obtained from Riemannian geometry by replacing the field $\R$ by
  $\C$ (e.g. a complex ``metric'' becomes a non-degenerate
  symmetric complex-bilinear form on the tangent space), and all
  classical results hold, naturally with the exception of those making
  use of partitions of unity. We will often omit the prefix
  ``complex-'', when referring to geometric objects, and we will always
  consider them, unless otherwise stated, in the framework of
  complex-Riemannian geometry.
\subsection{Conformal complex 4-manifolds}\label{ss.w21} Let {\bf M}
  be a 4-dimensional 
  complex manifold. A {\it conformal 
  structure} is defined, as in the real case \cite{gconf}, by a
  everywhere non-degenerate section $c$ of the
  complex bundle $S^2(T^*\mathbf{M})\otimes L^2$, where $L$ is a given line bundle
  of {\it scalars of weight 1}, and $L^4\simeq \kappa^{-1}$, the
  anti-canonical bundle of {\bf M}. (While on an oriented real manifold such a line
  bundle always exists, being topologically trivial, in the complex
  case the existence of $L^2$, a square root of the anti-canonical
  bundle, is submitted to some topological restrictions.) From now on,
  only holomorphic conformal structures will be considered, thus $L$
  is a holomorphic bundle and $c$ a holomorphic section of
  $S^2(T^*\mathbf{M})\otimes L^2$. (In fact, all we need to define the
  conformal structure $c$ on the 4-manifold {\bf M} is just the
  holomorphic bundle $L^2$; in odd dimensions the situation is
  different, see Section \ref{s.w7}.) 

As in the real case, $c$ is locally represented
  by symmetric bilinear forms on $T\mathbf{M}$, but global representative
  metrics do not exist, in general.

For each point $x\in \mathbf{M}$, there is an isotropy cone $C_x$ in the tangent
space $T_x\mathbf{M}$, who uniquely determines the conformal structure
$c$. In the associated projective space, $\P(T_x\mathbf{M})\simeq \cp3
$, 
the cone $C_x$ projects onto the non-degenerate quadratic surface
$\P(C_x)$, which 
is actually a ruled surface isomorphic to $\cp1\times\cp1$. We thus get 2
families of complex projective lines contained in $\P(C)$, that is, 2
families of isotropic 2-planes in $C\subset T\mathbf{M}$, respectively called
$\alpha$-planes, the other $\beta$-planes. This choice corresponds to
the choice of an ``orientation'' of {\bf M}. 
On a real 4-manifold an
orientation is chosen by picking a class of volume forms (which is not
possible in this complex framework) or by choosing one of the two
possible Hodge operators compatible with the conformal structure
$*:\Lambda^2\mathbf{M}\rightarrow \Lambda^2\mathbf{M}$ (which can also be done in our
complex case, \cite{p2}). As $*$ is a symmetric involution,
$\Lambda^2\mathbf{M}$ decomposes in $\Lambda^+{\bf
  M}\oplus\Lambda^-\mathbf{M}$ consisting in $\pm 1$-eigenvectors
of $*$, respectively called {\it self-dual} and {\it anti-self-dual}
2-forms; the isotropic vectors in $\Lambda^+{\bf M}$ and
$\Lambda^-{\bf M}$ are then exactly the
{\it decomposable}
elements $u\wedge v\in \Lambda^\pm\mathbf{M}$, with $u,v\in{\bf M}$.
\begin{defi} An $\alpha$-plane $F^\alpha$ (resp. a $\beta$-plane
  $F^\beta$ in $T\mathbf{M}$ is a 
  2-plane such that $\Lambda^2F^\alpha$ (resp. $\Lambda^2F^\beta$) is
  a is a self-dual (resp. anti-self-dual) isotropic line in
  $\Lambda^2\mathbf{M}$.
\end{defi}
\obs The $\alpha$- and $\beta$-planes can be interpreted in terms of
spinors. The structure group of the tangent bundle $T\mathbf{M}$ is restricted
to the conformal orthogonal complex group,
$CO(4,\C):=(O(4,\C)\times\C^*)/\{\pm\mathbf{1}\}$, where
$O(4,\C):=\{A\in GL(4,\C)|A^{\; t}\! A=\mbox{Id}\}$, by the given conformal
structure of {\bf M}. The choice of an orientation is the further
restriction of this group to the connected component of the identity,
$CO_0(4,\C):=SO(4,\C)\times\C^*$, where $SO(4,\C):=O(4,\C)\cap SL(4,\C)$.
Consider a local metric $g$ in the conformal class $ c$. We have then
locally defined $Spin$ structures, and associated $Spin$ bundles 
$V_+,V_-$, as in the real case \cite{ahs},\cite{sal}. They are rank 2
complex vector bundles, and for each local section of $L$ (i.e. a
metric in $c$), each of them is equipped with a (complex) symplectic
structure $\omega_+\in\Lambda^2V_+,\omega_-\in\Lambda^2V_-$,
respectively.  
Then we locally have $T\mathbf{M}\simeq V_+\otimes V_-$, and
$g=\omega_+\otimes\omega_-$, for the fixed metric $g\in c$.
$\alpha$-(resp. $\beta$-) planes are then nothing but the isotropic
2-planes obtained by fixing the first (resp. the second) factor in the
above tensor product:
\begin{prop}\cite{p2} An $\alpha$-plane, resp. $\beta$-plane $F\subset
  T_x\mathbf{M}$ is a complex 
  plane $\psi_+\otimes V_-$, resp $V_+\otimes \psi_-$, where
  $\psi_+\in V_+\smallsetminus\{0\}$, resp. $\psi_-\in
  V_-\smallsetminus\{0\}$. 
\end{prop}
The $\alpha$-planes in $T_x\mathbf{M}$ are, thus indexed by
$\P(V_+)_x$, and $\beta$-planes by $\P(V_-)_x$, and these
projective bundles are globally well-defined on {\bf M}, \cite{ahs}. 

\obs It is obvious that a change of orientation interchanges the
$\alpha$ and $\beta$-planes; the same is true for self-duality and
anti-self-duality, to be defined below.\medskip

For a local metric $g$ in $c$, we denote by $R^g$ its Riemannian
curvature, and by $W$ the Weyl tensor, i.e. the trace-free
component of $R^g$, which is known to be independent of the chosen
metric within the conformal class \cite{gconf}.
It splits into two components $W^+,W^-$, and the easiest way to see
that is the spinorial decomposition of the space of the curvature
tensors $\mathcal{R}\subset\Lambda^2 \otimes\Lambda^2$,
\cite{ahs},\cite{sal},\cite{st},  obtained from the relation
$T\mathbf{M}=V_+\otimes V_-$ and some of the Clebsch-Gordan identities \cite{sal}. 
$$\mathcal{R}=\mathcal{S}\oplus\mathcal{B}\oplus\mathcal{W}^+
\oplus\mathcal{W}^-,$$
where $\mathcal{S}$ is the complex line of scalar curvature tensors,
included in $\Lambda^2V_+\oplus\Lambda^2V_-$, $\mathcal{B}=S^2V_+\otimes
S^2V_-$ is the space of trace-free Ricci tensors, and $\mathcal{W}^+=S^4V_+$,
$\mathcal{W}^-=S^4V_-$ are the spaces of self-dual,
resp. anti-self-dual Weyl tensors (where $S^pV_\pm$ denotes the
$p$-symmetric power of $V_\pm$). 

The curvature $R^g$ restricted to 
  any $\alpha$-plane $F$ yields a weighted bilinear symmetric form
  $R^F$ on $\Lambda^2F$, i.e. a section in $L^2\otimes(\Lambda^2
  F\otimes\Lambda^2 F)^*$:
$$(g,X\wedge Y)\longmapsto g(R^g(X,Y)X,Y).$$
\begin{prop}
  The (weighted) bilinear form $R^F$ depends only on the self-dual Weyl
  tensor, and this one is completely determined by the (weighted)
  values of $R^F$ for all $\alpha$-planes $F$. 
\end{prop}

We have the same result for $\beta$-planes.

\begin{proof} Let $F=V_+\otimes\psi_-$ be an $\alpha$-plane, and let
$X=\psi_+\otimes\fy_1, Y=\psi_+\otimes\fy_2\in F$, and suppose, for
simplicity, that $\omega_-(\fy_1,\fy_2)=1$, so $X\wedge Y\in\Lambda^2F$
is identified to the element $\psi_+\otimes\psi_+\in S^2V_+$. Then it is easy 
to see that $R^F$, evaluated on $X\wedge Y$, is 
nothing but the evaluation of $R\in S^2(\Lambda^2\mathbf{M})\supset\mathcal{R}$
on $(X\wedge Y)\otimes(X\wedge Y)\simeq \psi_+\otimes\psi_+\otimes
\psi_+\otimes\psi_+\in S^4V_+$, which depends only on the positive (or
self-dual) part of the Weyl tensor. 
To prove the second assertion, we remark that $W^+$, being a 
quadrilinear symmetric form on $V_+$, can be identified with a polynomial
of degree 4 on $V_+$, which is determined by its values.
\end{proof}

\begin{defi} A conformal structure $c$ on a 4-manifold {\bf M} is called 
self-dual (resp. anti-self-dual) iff $W^-=0$ (resp. $W^+=0$).
\end{defi}

\obs In general, geodesics on a conformal manifold depend on the
chosen metric, with the exception of the isotropic ones (or {\it
  null-geodesics}). Therefore the existence of totally geodesic
surfaces tangent to $\alpha$- (resp. $\beta$-) planes is a property of
the conformal structure alone. 
\subsection{Twistor spaces}\label{ss.w22}
\begin{defi}
  An $\alpha$-surface (resp. $\beta$-surface) $\alpha\subset \mathbf{M}$ is a maximal, totally geodesic
  surface in {\bf M}, whose tangent space in any point is an $\alpha$-plane
  (resp. $\beta$-plane).
\end{defi}
On the other hand, any totally geodesic, isotropic surface in {\bf M}
is included in an $\alpha$- or in a $\beta$-surface.
\begin{defi}\cite{p1},\cite{p2} If, in any point $x\in \mathbf{M}$, and for any
  $\alpha$- 
  (resp. $\beta$-) plane $F\subset T_x\mathbf{M}$, there is a $\alpha$-
  (resp. $\beta$-) surface tangent to $F$ at $x$, we say that the
  family of $\alpha$- (resp. $\beta$-) planes is integrable. 
\end{defi}

%
%
%
%
%
%

\begin{thrm.} \cite{ahs},\cite{p2}
The family of $\alpha$- (resp. $\beta$-) planes of a conformal
4-manifold $(\mathbf{M},c)$  integrable if and only if the conformal structure
$c$ is anti-self-dual (resp. self-dual). 
\end{thrm.}

The integrability of $\alpha$-planes is equivalent to the integrability
(in the sense of Frobenius) of a distribution $H^\alpha$ of 2-planes on the
total space of the projective bundle $\P(V_+)$. More precisely, let $g$
be a local metric in the conformal class $c$, and let $\nabla$ be its 
Levi-Civita connection. $\nabla$ induces a connection in the bundle  
$\P(V_+)$, thus a horizontal distribution $H$, isomorphic to $T\mathbf{M}$
via the bundle projection. Let $H^\alpha$ be the 2-dimensional
subspace of $H_F$ --- where $F\in \P(V_+)$ is an $\alpha$-plane in
$T_x\mathbf{M}$ --- which projects onto $F\subset T_x\mathbf{M}$. It can 
be easily shown (as in \cite{p2}, see also \cite{ahs}) that the  
``tautological'' 2-plane distribution $H^\alpha$ is independent of the metric $g$.
Then $\alpha$-surfaces are canonically lifted as integrable manifolds
of the distribution $H^\alpha$. For a geodesically convex open set of
{\bf M}, one can prove (see \cite{leb5}) that the space of these integrable
leaves is a complex 3-manifold. (This point of view is closely related to
the one of \cite{ahs}, about the integrability of the canonical
almost complex structure of the real twistor space.)

The same remark can be made about $\beta$-surfaces.

\obs The existence, for any point $x\in \mathbf{M}$, of an $\alpha$-surface containing $x$
does not imply, in general, the integrability of the family of $\alpha$-planes :
in the conformal self-dual (but not anti-self-dual) manifold
$\mathbf{M}=\cp2\times(\cp2)^*\smallsetminus\mathcal{F}$ (the complexification of
$\cp2$, \cite{ahs}), the surfaces $(\{x\}\times(\cp2)^*)\cap \mathbf{M}$ and
$(\cp2\times\{y\})\cap \mathbf{M}$ are all $\alpha$-surfaces, see Section \ref{s.w8}.

\obs In the real framework, the twistor space of a real Riemannian
4-manifold $M^\R$ is the total space $Z^\R$ of the $S^2$-bundle of
almost-complex structures on $TM^\R$, compatible with the conformal
structure and the (opposite) orientation; it admits a natural
almost-complex structure $\mathcal{J}$, equal, in $J\in Z^\R$, to the complex
structure of the fibers on the {\it vertical} space $T^\vee_J Z^\R$,
and to $J$ itself on the horizontal space (induced by the Levi-Civita
connection). Such a complex structure $J$ is equivalent to an
isotropic complex 2-plane in $TM\otimes\C$, thus to an $\alpha$- or
$\beta$-surface (depending on the conventions), which becomes then the space of
vectors of type $(1,0)$ for $J$; as the integrability of the
almost-complex structure $\mathcal{J}$ can be expressed as the
Frobenius condition applied to $T^{(1,0)}Z^\R$, it is equivalent to
the integrability of the family of $\alpha$-, resp. $\beta$-planes.
\medskip

The {\it Penrose construction} associates to an (anti-)self-dual
manifold {\bf M} the space $Z$ of $\alpha$- (resp. $\beta$-)surfaces of {\bf M};
we have seen above that $Z$ admits complex-analytic maps, but it may
be non-Hausdorff. This is why we need to introduce the following
condition , see also \cite{leb5}:
\begin{defi}
  An (anti-) self-dual manifold {\bf M} is called {\rm civilized} iff the
  space $Z^\alpha$ (resp. $Z^\beta$) of integral leaves of the
  distribution $H^\alpha$  (resp. $H^\beta$) in $\P(V_+)$
  (resp. $\P(V_-)$) is a complex   3-manifold, and the projection
  $p^+:\P(V_+)\rightarrow Z^\alpha$ (resp. $p^-: \P(V_-)\rightarrow
  Z^\beta$) is a submersion.

In this case, the manifold $Z^\alpha$ (resp. $Z^\beta$) --- which is
the space of $\alpha$-surfaces (resp. $\beta$-surface) of {\bf M} --- is called the $\alpha$-
(resp. the $\beta$-)twistor space of {\bf M}.
\end{defi}
\bigskip

From now on, we suppose that $(\mathbf{M},c)$ is a self-dual complex analytic
4-manifold. As any point $x\in \mathbf{M}$ has a geodesically convex
neighborhood $U$ \cite{ehr} (which is, therefore, civilized), we can
construct $Z^U$, the $\beta$-twistor space (for short,
twistor space) of $U$. For the infinitesimal results of this paper
(from Sections \ref{s.w3},\ref{s.w4},\ref{s.w6} and \ref{s.w7}), we
will suppose (with no loss of 
generality) that {\bf M} is civilized (for example, by replacing {\bf M} by $U$).
\bigskip

We recall now the correspondence between differential
geometric objects on {\bf M} and complex analytic objects on its twistor
space, $Z$, \cite{ahs},\cite{p2}, see also
\cite{leb1},\cite{leb3},\cite{sal}. 

$\beta$-surfaces $\beta\subset \mathbf{M}$ correspond to points $\bar\beta\in Z$, by
definition, and the set of $\beta$-surfaces passing through a point
$x\in \mathbf{M}$ is a 
complex projective line $Z_x$,
with normal bundle isomorphic (non-canonically) to
$\mathcal{O}(1) \oplus\mathcal{O}(1)$ (where $\o1$ is the dual of the
{\it tautological} bundle $\o{-1}$ on $\cp1$) \cite{ahs},\cite{p2}, see also
\cite{besse}.

In fact, this family of complex projective lines in $Z$ permits us to
recover {\bf M} and its conformal structure, at least locally, by the {\it
reverse Penrose construction}: 
The normal bundle $N_x$ of a line $Z_x$
in $Z$ has the property $H^1(N_x,\mathcal{O})=0$, thus, by a theorem
of Kodaira \cite{kod}, the 
space $\mathbf{M}_0$ of projective lines in $Z$ having the above normal bundle is a
smooth complex manifold, whose tangent space at a point $x\simeq
Z_x\subset Z$ is canonically isomorphic to the space of global
sections of the normal bundle $N_x$ of $Z_x$ (thus  $\mathbf{M}_0$ has dimension 4).
The conformal structure of $\mathbf{M}_0$ is described by its tangent cone,
which corresponds to the sections of $N_x$ having at least one zero
(as such a section decomposes as 2 sections of $\mathcal{O}(1)$, the
vanishing condition means that they both vanish at the same point, which
is a quadratic condition on the sections of $N_x$). We thus get a
conformal diffeomorphism from {\bf M} to an open set of $\mathbf{M}_0$.

\subsection{Ambitwistor spaces}\label{ss.w23}
We remark that $\P(V_-)$ is an open set of the projective tangent
bundle of $Z$, as $Z$ is the space of leaves of $\P(V^-)$, but it is
important to note that, 
in general, the reverse inclusion is not true (i.e. not any direction
in $Z$ is tangent to a line corresponding to a point in {\bf M}, or,
equivalently, $\beta$-surfaces are not compact $\cp2$'s, in general, see Section
\ref{s.w5}). 

For example, if $\mathbf{M}=\m$ (with the notations in Section \ref{s.w8}),
$\P(V_-)$ is an open subset in the $\cp2$-bundle $\P(TZ)\rightarrow
Z$, consisting in the set of directions transverse to the {\it contact
 structure} of $Z$ (see Section \ref{ss.w84}). $\P(V_-)$ is, thus, in this
case, a rank 2 affine bundle over $Z$. 

Another canonical $\cp2$-bundle on $Z$, that is
$\P(T^*Z)\rightarrow Z$, leads to the {\it ambitwistor space} $B$, by definition
the space of null-geodesics of {\bf M} \cite{leb3}. It is an open
set of the projective cotangent bundle of $Z$ (or, equivalently, the
Grassmannian of 2-planes in $TZ$) \cite{leb3} (more precisely, a plane
$F\subset T_{\bar\beta}Z$ 
corresponds to a null-geodesic $\gamma\subset \mathbf{M}$ (contained
in $\beta$) if it is 
tangent to at least one projective line $Z_x$, corresponding to a
point $x\in \mathbf{M}$).  

To see that, let $x$ be a point in
$\mathbf{M}$, $\beta$ a $\beta$-surface passing through $x$, i.e. $\bar\beta\in Z$ and
$Z_x$ contains $\bar\beta$; let $F\subset T_{\bar\beta} Z$ be a plane tangent to
$Z_x$. As small deformations of $Z_x$ still correspond to points of
{\bf M}, we consider those rational curves which 
are tangent to $F$. They correspond to a (continuous) set of points on a
curve $\gamma\subset \beta$, that will turn out to be a null-geodesic. Indeed, all
we have to prove is $\ddot\gamma=0(\mbox{mod}\;\dot\gamma)$, and
$\dot\gamma_x$ corresponds to a section $\eta$ of $N_x$, vanishing at
$\bar\beta\in Z_x$; as $N_x\simeq\o1$, $\eta$ is determined by its
derivative at $\bar\beta$, which is a linear map
$T_{\bar\beta}\rightarrow F/T_{\bar\beta}$ (the infinitesimal
deformation of the direction of $Z_x$ within $F$). As the points of
$\gamma$ correspond to lines tangent to $F$, we have that
$\ddot\gamma_x$ corresponds to a section of $N_x$ collinear to $\eta$,
thus $\gamma$ verifies the equation of a (non-parameterized)
geodesic. See \cite{leb1}, \cite{leb5}, and Section \ref{s.w4} for
details.
\smallskip

{\noindent\bf Example. } The space of null-geodesics of $\mathbf{M}=\M$ is the total
space of a $\C\times\cp1$-bundle over $Z=\mathcal{F}$, the flag
manifold (see Section \ref{s.w8}); a 2-plane $F\subset T_{(L,l)}\mathcal{F}$
which corresponds to a null-geodesic in {\bf M} is identified either to a projective
diffeomorphism $\varphi:\P(l)\rightarrow\P(L^o)$ (Section \ref{ss.w84}, case
3), or to a point $A\subset l,\, A\neq L$, resp. a plane $a$
containing $L$, and different from $l$ (Section \ref{ss.w84}, cases
{\bf 2} and ${\bf 2'}$). 

\section[The field of $\alpha$-cones]{The structure of the ambitwistor
  space and the field of $\alpha$-cones}\label{s.w3} 
{\bf Conventions.} Except for some results in Section \ref{s.w5}, we will
consider {\bf M} to be a self-dual civilized 4-manifold, i.e. the
(twistor) space $Z$ of $\beta$-surfaces of {\bf M} is a Hausdorff smooth complex
3-manifold, and the projection $\P(V_-)\rightarrow Z$ is a submersion
(e.g. {\bf M} is geodesically convex), see \cite{leb5}.

We will frequently identify, following the deformation theory of
Kodaira (see \cite{kod}), the vectors in $T_x\mathbf{M}$ with sections in the
normal bundle $N(Z_x)$ of the projective line $Z_x$ in $Z$.
 
We also consider the space of null-geodesics $B$, as an open subsetset
of $\P(T^*Z)$.

For a null-geodesic $\gamma$, resp. a $\beta$-surface $\beta\subset
\mathbf{M}$, we denote by 
$\bar\gamma$, resp. $\bar\beta$, the corresponding point in $B$,
resp. $Z$.
 
\subsection{$\alpha$- and $\beta$-cones on the ambitwistor space}\label{ss.w31}

The vectors on $B$ can be expressed in terms
of infinitesimal deformations of geodesics of
{\bf M} (Jacobi fields). More precisely:
$$T_{\bar\gamma} B\simeq\j^\perp_\gamma/\j^\gamma_\gamma,$$
where, for a null-geodesic $\gamma$, $\j_\gamma^\perp$ is the space of Jacobi
fields $J$ such that $\nabla_{\dot\gamma}J\perp\dot\gamma$, and
$\j^\gamma_\gamma$ is its subspace of Jacobi fields ``along''
$\gamma$, i.e. $J\in\C\dot\gamma$ in any point of the geodesic.

\obs A class in $\j^\perp_\gamma/\j^\gamma_\gamma$ is represented by
Jacobi fields yielding the same local section of the
normal bundle $N(\gamma)$ of $\gamma$ in {\bf M}. This is equivalent to
the following obvious fact: 

\begin{lem}
  The kernel of the natural application $\j^\perp_\gamma\rightarrow
  N(\gamma)$ is $\j^\gamma_\gamma$.
\end{lem}

As a consequence, Jacobi fields on $\gamma$ induce particular local
sections in $N(\gamma)$, which turn out to be solutions of a
differential operator of order 1 on $N(\gamma)$, see Section \ref{s.w4}. 
\medskip

The conformal geometry of {\bf M} induces a particular structure on $B$:
we describe it in order to obtain an expression of $W^+$ in terms of
the geometry of the (ambi-\nolinebreak{})\linebreak[0]{}twistor space.
\bigskip

We have a canonical hyperplane $V_{\bar\gamma}$ in $T_{\bar\gamma} B$,
defined by 
$$ V-{\bar\gamma}:=\j^{\perp\perp}_\gamma/\j^\gamma_\gamma,$$
where $\j^{\perp\perp}_\gamma$ is the set of Jacobi fields $J$ everywhere
orthogonal to $\dot\gamma$ (i.e. $\nb_{\dot\gamma}J\perp\dot\gamma$
{\it and} $J\perp\dot\gamma$). 

We deine now two fiels of cones in $TB$, both contained in
$V-{\bar\gamma}$:
\begin{defi} Let $\gamma$ be a null-geodesic in {\bf M}, and, for each
  point $x\in\gamma$, let $F^\beta_x$ be the $\beta$-plane containing
  $\dot\gamma_x$. The (infinitesimal) $\beta$-cone 
  $V^\beta_{\bar\gamma}$ at $\bar\gamma\in B$ is defined as follows:
$$  V^\beta_{\bar\gamma}:=\j^\beta_\gamma/\j^\gamma_\gamma\subset
\j^{\perp\perp}_\gamma/\j^\gamma_\gamma\simeq V_{\bar\gamma} \subset
T_{\bar\gamma}B,$$
where $\j^\beta_\gamma$ is the set of Jacobi fields $J$ on $\gamma$
satisfying the condition
$$\exists\, x\in\gamma \mbox{ such that }J_x=0 \mbox{ and }
(\nb_{\dot\gamma}J)_x\in F^\beta_x.$$
\end{defi}

\begin{prop}
The $\beta$-cone $V^\beta_{\bar\gamma}$ is {\em flat}, i.e. it is included in 
the 2-plane $F_{\bar\gamma}^\beta$ consisting of
Jacobi fields contained in the $\beta$-plane defined by $\dot\gamma$ in each
point of it.
\end{prop}
\begin{proof}
We have to prove that $\j^\beta_\gamma$ is  included in
$\bar\j^\beta_\gamma$, defined as follows:
$$\bar\j^\beta_\gamma:=\{J\mbox{ Jacobi field on }\gamma\; |\; J_x,\dot
J_x\in F^\beta_x,\, \forall x\in\gamma\}.$$
We will prove that $\j^\beta_\gamma\subset\bar\j^\beta_\gamma$,
therefore it will follow that the latter is non-empty, and is a linear
space of dimension 2.

We denote by $J^0$ the parallel
displacement, along $\gamma$,  of a non-zero vector in 
$F^\beta_x$, transverse to $\dot\gamma$. Then $J^0\in
T\beta|_\gamma\smallsetminus T\gamma$, because $\gamma$ is included in
the is totally geodesic surface $\beta$,
thus we can characterize $F^\beta_y$ as the set $\{X\in T_yM\; |\; X\perp\dot\gamma,
X\perp J^0\}$, for any $y\in\gamma$. We then observe that 
$$\dot\gamma .\langle \dot
J,J^0\rangle =\langle R(\dot\gamma,J)\dot\gamma, J^0\rangle =\langle R(\dot\gamma,J^0)\dot\gamma,
J\rangle =k\langle J^0,J\rangle ,$$ 
because $R(\dot\gamma,J^0)\dot\gamma$ is in $F^\beta$, thus
$R(\dot\gamma, J^0)\dot\gamma=h\dot\gamma +kJ^0.$ So the scalar
function $\langle J,J^0\rangle $ satisfies to a linear second order equation, hence
it it determined by its initial value and derivative. It follows then
that it is identically zero, thus $J\in F^\beta$ everywhere, as
claimed. \end{proof}

Another subset in $T_{\bar\gamma} B$ is the $\alpha$-cone $V^\alpha_{\bar\gamma}$,
defined as follows: 
\begin{defi} Let $\gamma$ be a null-geodesic in {\bf M}, and, for each
  point $x\in\gamma$, let $F^\alpha_x$ be the $\alpha$-plane containing
  $\dot\gamma_x$. The (infinitesimal) $\alpha$-cone 
  $V^\alpha_{\bar\gamma}$ at $\bar\gamma\in B$ is defined as follows:
$$  V^\alpha_{\bar\gamma}:=\j^\alpha_\gamma/\j^\gamma_\gamma\subset
\j^{\perp\perp}_\gamma/\j^\gamma_\gamma\simeq V_{\bar\gamma} \subset
T_{\bar\gamma}B,$$
where $\j^\alpha_\gamma$ is the set of Jacobi fields $J$ on $\gamma$
satisfying the condition
$$\exists x\in\gamma \mbox{ such that }J_x=0 \mbox{ and }
(\nb_{\dot\gamma}J)_x\in F^\alpha_x.$$
\end{defi}

It is important to note that, in general, the
projective curves $\P(V^\alpha_{\bar\gamma})$ and
$\P(V^\beta_{\bar\gamma})$ are non compact, as each of them 
corresponds to the set of points on $\gamma$, which is non-compact, in
general. The field of $\alpha$-cones on $B$ is the object of main
interest in this paper. We may already guess that its flatness
(i.e. the situation when $V^\alpha_{\bar\gamma}$ is a subset in a 2-plane)
can be related to some special property of the conformal structure of
{\bf M}.

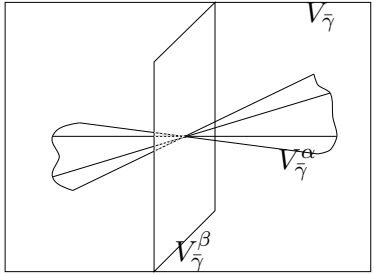
\begin{center}\begin{picture}(0,0)%
\epsfig{file=weyl/cone.pstex}%
\end{picture}%
\setlength{\unitlength}{0.00015700in}%
\begingroup\makeatletter\ifx\SetFigFont\undefined
\def\x#1#2#3#4#5#6#7\relax{\def\x{#1#2#3#4#5#6}}%
\expandafter\x\fmtname xxxxxx\relax \def\y{splain}%
\ifx\x\y   
\gdef\SetFigFont#1#2#3{%
  \ifnum #1<17\tiny\else \ifnum #1<20\small\else
  \ifnum #1<24\normalsize\else \ifnum #1<29\large\else
  \ifnum #1<34\Large\else \ifnum #1<41\LARGE\else
     \huge\fi\fi\fi\fi\fi\fi
  \csname #3\endcsname}%
\else
\gdef\SetFigFont#1#2#3{\begingroup
  \count@#1\relax \ifnum 25<\count@\count@25\fi
  \def\x{\endgroup\@setsize\SetFigFont{#2pt}}%
  \expandafter\x
    \csname \romannumeral\the\count@ pt\expandafter\endcsname
    \csname @\romannumeral\the\count@ pt\endcsname
  \csname #3\endcsname}%
\fi
\fi\endgroup
\begin{picture}(12174,9291)(214,-8398)
\put(9631,-4831){\makebox(0,0)[lb]{\smash{\SetFigFont{11}{13.2}{rm}
$V^\alpha_{\bar\gamma}$
}}}
\put(10576,209){\makebox(0,0)[lb]{\smash{\SetFigFont{11}{13.2}{rm}
$V_{\bar\gamma}$
}}}
\put(6076,-7981){\makebox(0,0)[lb]{\smash{\SetFigFont{11}{13.2}{rm}
$V^\beta_{\bar\gamma}$
}}}
\end{picture}
\end{center}

\obs We have seen that $V^\beta_{\bar\gamma}$ is included in
the 2-plane $F^\beta_{\bar\gamma}$, i.e. the condition  $J_x=0,\dot
J_x\in F_x$ can be generalised to the linear condition $J,\dot J\in
F^\beta$, but there is no canonical way
of supplying the ``missing'' points of $\gamma$ with some appropriate
Jacobi fields in order to ``complete'' $V^\alpha_{\bar\gamma}$ as in the
$\beta$-cones case. This would be possible, for example, if
$\P(V^\alpha_{\bar\gamma})$ would be an open subset in a 
projective line. But the defect of $V^\alpha_{\bar\gamma}$ to be part of a
2-plane is measured by its {\it projective curvature}, and we will see
in Section \ref{s.w4} that the vanishing of the latter implies the vanishing of
$W^+$ (Theorem \ref{t1}).
$B^\alpha_{\bar\gamma,x}$
\subsection{Integral $\alpha$-cones in $Z$ and $B$}\label{ss.w32}
We study now the field of $\alpha$-cones of $B$ in relation with $Z$ and the
canonical projection $\pi:B\rightarrow Z$. First, we note that there
are complex projective lines in $B$ tangent to the
directions in $V^\alpha_{\bar\gamma}$: 
\begin{defi} Let $\bar\gamma\in B$, $x\in\gamma$ a point on the
  null-geodesic $\gamma$; let $F^\alpha_x$ be the $\alpha$-plane
  tangent to $T_x\gamma$. The rational curve
  $B^\alpha_{\bar\gamma,x}$ in $B$ (containing $\bar\gamma$), is
  by definition the set of null-geodesics passing through $x$ and
  tangent to $F^\alpha_x$.
\end{defi}

The curves $B^\alpha_{\bar\gamma,x}, \, x\in\gamma$ are projected by
$\pi$ onto the complex lines $Z_x$
through $\bar\beta$  (corresponding to the $\beta$-surface $\beta$
containing $\gamma$) tangent to the
2-plane $F^\gamma$.

On the other hand, it is easy to see that the complex projective lines
$B^\beta_{\bar\gamma,x}$ (defined in an analogous way to $B^\alpha_{\bar\gamma,x}$), which are tangent to (an open set of the directions of)
$V^\beta_{\bar\gamma}$, are contained in the fibers of $\pi$. In fact, they
coincide with some of the projective lines passing through the point
$\gamma\in\P(T^*_{\bar\beta} Z)\simeq \cp2$.
\begin{defi}
The {\it integral $\alpha$-cones} in $B$, resp. $Z$ are defined by: 
$$B^\alpha_{\bar\gamma}:=\bigcup_{x\in\gamma}B^\alpha_{\bar\gamma,x}\
\mbox{($\beta$-cone in $B$)};\
Z^\gamma:=\bigcup_{x\in\gamma}Z_x\ \mbox{($\beta$-cone in $Z$)}.$$
\end{defi}
We intend to prove that $B^\alpha_{\bar\gamma}$ is the {\it canonical lift}
of $Z^\gamma$ (see Proposition \ref{p4}).
We know that $\pi( B^\alpha_{\bar\gamma})=Z^\gamma$. We have then the
following:

\begin{prop}\label{p1}
 Except for the vertices $\bar\gamma\in B^\alpha_{\bar\gamma}$ and
  $\bar\beta\in Z^\gamma$, the two integral cones $B^\alpha_{\bar\gamma}$
  and $Z^\gamma$ are smooth, immersed, surfaces of $B$, resp. $Z$.
\end{prop}
\begin{proof}
The open set of $B$ which is the space of null-geodesics of {\bf M} can be viewed as
the space of integral curves of the {\it geodesic distribution} $G$ of
lines in $\P(C)$, the total space of the fibre bundle of isotropic
directions in $T\mathbf{M}$. $G_v$ is defined as the horizontal lift (for the
Levi-Civita connection on {\bf M}) of $v$, which is an isotropic line in
$T_x\mathbf{M}$. This definition is independent of the chosen metric
and connection \cite{leb5}, and, by integrating this distribution (as
{\bf M} is civilized), we get a
holomorphic map $p:\P(C)\rightarrow B$, where (an open set of) $B$ is
the space of leaves of this foliation. This map can be used to
compute the normal bundle of $B^\alpha_{\bar\gamma,x}$, $N(B^\alpha_{
\bar\gamma,x})$, see \cite{leb1},\cite{leb3},\cite{leb5}. 

Indeed, we have lines $C^\alpha_{\gamma,x}\in \P(C)_x$, such that
$\dot\gamma_x\in C^\alpha_{\gamma,x}$, 
which project onto $B^\alpha_{\bar\gamma,x}$, thus we get the following
exact sequence of normal bundles:
$$0\rightarrow N(C^\alpha_{\gamma,x};p^{-1}(B^\alpha_{\bar\gamma,x}))
\rightarrow N(C^\alpha_{\gamma,x};\P(C))\rightarrow
N(B^\alpha_{\bar\gamma,x};B)\rightarrow 0,$$
where we have written the ambient spaces of the normal bundles on the
second position. The central bundle is trivial ($C^\alpha_{\gamma,x}$
is trivially embedded in $\P(C)_x\simeq\cp1\times\cp1$, which is
trivially embedded in $\P(C)$ as a fibre), and it is easy to check
that the left hand bundle is isomorphic to the tautological bundle over
$\cp1$, $\o{-1}$. This proves that $N(B^\alpha_{\bar\gamma,x};B)\simeq
\o0\oplus\o0\oplus\o1$, in particular the conditions in the
completeness theorem of Kodaira \cite{kod} are satisfied. Thus the
lines in the integral $\alpha$-cone $B^\alpha_{\bar\gamma}$ form an analytic
subfamily of the family $\{B^\alpha_{\bar\gamma,x}\}_{\bar\gamma\in
  B,x\in\gamma\subset M}$, that correspond to the sections of the
normal bundle of $B^\alpha_{\bar\gamma,x}$, vanishing at $\bar\gamma\in
B$, or, equivalently, to the points $x$ of $\gamma\subset \mathbf{M}$. 
\smallskip

But, in order to prove the smoothness of $B^\alpha_{\bar\gamma}\smallsetminus\{
\bar\gamma\}$,  we first remark that the surface
$C^\alpha_\gamma\subset \P(C)$, defined as follows, is smooth:
$$C^\alpha_\gamma:= \{v\in \P(C)_x|x\in\gamma,\; v\subset
F^\alpha_\gamma\},$$
where $F^\alpha_\gamma$ is the $\alpha$-plane containing
$\dot\gamma$. $C^\alpha_\gamma$ is smooth, and
$p(C^\alpha_\gamma)=B^\alpha_{\bar\gamma}$. We note now that
$C^\alpha_\gamma$ is everywhere, with the exception of the points of
$p^{-1}(\bar\gamma)$, transverse to the fibers of the submersion
$p:\P(C)\rightarrow B$. We may conclude that
$B^\alpha_{\bar\gamma}\smallsetminus\{\gamma\}$ is a smooth analytic submanifold
of $B$ (not closed).

We can use similar methods to prove that
$Z^\gamma\smallsetminus\{\bar\beta\}$ is an immersed submanifold of $Z$ (by
using the projection $\pi:B\rightarrow Z$). 
\end{proof}

There is another argument for this latter claim, 
which gives the tangent
space to $Z^\gamma$ in any point:

We see $Z^\gamma$ as the ``trajectory'' of a 1-parameter deformation of
$Z_x$: we fix $\bar\beta$ and we ``turn'' $Z_x$ around $\bar\beta$ by
keeping it tangent to $F^\gamma$. The trajectory of this deformation is smooth in
$\zeta\in Z^\gamma\smallsetminus\beta$ iff any non-identically zero
section $\nu$ of the normal bundle 
$N(Z_x)$ corresponding to this 1-parameter deformation) does not
vanish at $\zeta$. In particular, the tangent space $T_\zeta Z^\gamma$
is spanned by $T_\zeta Z_x$ and $\nu(\zeta)$.

But the sections $\nu$ generating this deformation are the sections of
$N(Z_x)$ vanishing at $\bar\beta$, and they vanish at only one point (and
even there, only at order 0) unless they are identically zero, because
$N(Z_x)\simeq \o1\oplus\o1$. 

\obs The values of these sections in the
points of $Z_x$ other than $\bar\beta$, plus their derivatives in $\bar\beta$
(well-defined as they all vanish at $\bar\beta$), define a 1-dimensional
subbundle of $N(Z_x)$, which is isomorphic to $\o1$. In fact, we have
a 1-1 correspondence between the subbundles of $N(Z_x)$ isomorphic to
$\o1$ and the 2-planes in $T_{\bar\beta} Z$. Then, the space of holomorphic
sections of such a bundle is a linear space of dimension 2, consisting
in a family of sections of $N(Z_x)$ vanishing on {\it different}
points of $Z_x$. Thus we get a 2-plane
$F^\alpha$ of isotropic vectors in $T_x\mathbf{M}$, which is easily
seen to be an $\alpha$-plane, as the 
$\beta$-plane $F^\beta_x=T_x\beta$ consists in the set of all sections
of $N(Z_x)$ vanishing at 
$\bar\beta$ (we have $F^\alpha_x\cap T_x\beta =T_x\gamma$). The
tangent space to $Z^\gamma$, in a point $\zeta\in Z_x$, is spanned by
the subbundle of $N(Z_x)$ (isomorphic to $\o1$ --- see above), defined
by the isotropic vectors $v\in F^\alpha_x$. If $\gamma^\zeta$ is the
null-geodesic generated 
by $v^\zeta$, we conclude that $T_\zeta Z^\gamma$ is the 2-plane determined by
$\gamma^\zeta$, and that $\zeta=\pi(\overline{\gamma^\zeta})$.

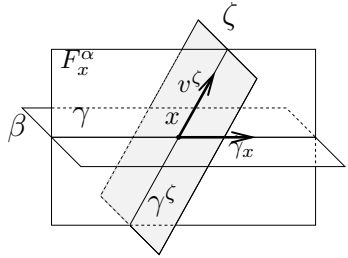
\begin{center}\input{weyl/zeta.pstex_t}\end{center}

{\noindent\bf Example. } If $\mathbf{M}=\M$, then the integral $\alpha$-cone
$Z^\gamma$ in $Z$, for $\gamma\equiv F^\gamma=F^\varphi\subset
T_{(L,l)}Z$ (where $\varphi:\P(l)\rightarrow\P(L^o)$ is a projective
diffeomorphism), is the (smooth away from the vertex $(L,l)$) surface
$\{(S,\varphi(s))|S\neq L, s\neq l, 
\varphi (s\cap l)=S\}$. Its compactification (by adding the {\it
  special cycle} $\bar Z_{(L,l)}$) is singular (Section \ref{ss.w84}).

As any smooth surface in $Z$ has a canonical lift in $B=\P(T^*Z)$, we
get:

\begin{prop}\label{p4} The integral $\alpha$-cone $B^\alpha_{\bar\gamma}$ is the
  canonical lift of the integral $\alpha$-cone on $Z$, $Z^\gamma$.
\end{prop}

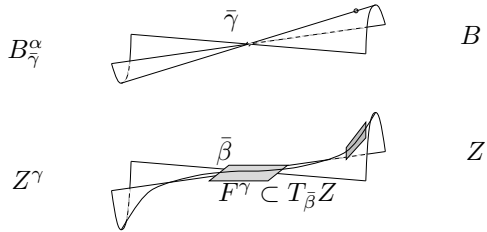
\begin{center}\begin{picture}(0,0)%
\epsfig{file=weyl/acones.pstex}%
\end{picture}%
\setlength{\unitlength}{0.00022700in}%
\begingroup\makeatletter\ifx\SetFigFont\undefined
\def\x#1#2#3#4#5#6#7\relax{\def\x{#1#2#3#4#5#6}}%
\expandafter\x\fmtname xxxxxx\relax \def\y{splain}%
\ifx\x\y   
\gdef\SetFigFont#1#2#3{%
  \ifnum #1<17\tiny\else \ifnum #1<20\small\else
  \ifnum #1<24\normalsize\else \ifnum #1<29\large\else
  \ifnum #1<34\Large\else \ifnum #1<41\LARGE\else
     \huge\fi\fi\fi\fi\fi\fi
  \csname #3\endcsname}%
\else
\gdef\SetFigFont#1#2#3{\begingroup
  \count@#1\relax \ifnum 25<\count@\count@25\fi
  \def\x{\endgroup\@setsize\SetFigFont{#2pt}}%
  \expandafter\x
    \csname \romannumeral\the\count@ pt\expandafter\endcsname
    \csname @\romannumeral\the\count@ pt\endcsname
  \csname #3\endcsname}%
\fi
\fi\endgroup
\begin{picture}(10485,5197)(361,-5246)
\put(5311,-556){\makebox(0,0)[lb]{\smash{\SetFigFont{10}{12.0}{rm}
$\bar\gamma$
}}}
\put(5131,-3616){\makebox(0,0)[lb]{\smash{\SetFigFont{10}{12.0}{rm}
$\bar\beta$
}}}
\put(5176,-4561){\makebox(0,0)[lb]{\smash{\SetFigFont{10}{12.0}{rm}
$F^\gamma\subset T_{\bar\beta} Z$
}}}
\put(10846,-3616){\makebox(0,0)[lb]{\smash{\SetFigFont{10}{12.0}{rm}
$Z$
}}}
\put(10711,-961){\makebox(0,0)[lb]{\smash{\SetFigFont{10}{12.0}{rm}
$B$
}}}
\put(361,-1276){\makebox(0,0)[lb]{\smash{\SetFigFont{10}{12.0}{rm}
$B^\alpha_{\bar\gamma}$
}}}
\put(451,-4291){\makebox(0,0)[lb]{\smash{\SetFigFont{10}{12.0}{rm}
$Z^\gamma$
}}}
\end{picture}
\end{center}

\obs Basically, this lift can only be defined for
$Z^\gamma\smallsetminus\{\bar\beta\}$, but, in this special case, it can be
extended by continuity to $\bar\beta$. Of course, the smoothness of the
lifted surface can only be deduced away from the vertex $\bar\gamma$ (from
the smoothness of $Z^\gamma\smallsetminus\{\bar\beta\}$).

\section[The projective curvature of $V^\alpha_\gamma$ and the Weyl tensor of {\bf M}]{The projective curvature of the $\alpha$-cone $V^\alpha_\gamma$ and the self-dual Weyl tensor $W^+$ on {\bf M}}\label{s.w4}

As noted in Section \ref{s.w3}, we intend to find a relation
 between the ``curvature'' of the $\alpha$-cone
$V^\alpha_{\bar\gamma}$ (its non-flatness) and the Weyl tensor $W^+$ of
$(\mathbf{M},c)$. We begin by defining the {\it projective curvature } of
$V^\alpha_\gamma$: A projective structure on a manifold $X$ is an
 equivalence class of linear connections yielding the same geodesics. 
In such a space, we
can define the {\it projective curvature} of a curve $S$ in a point
$\sigma$ as the linear application $k:T_\sigma S\otimes T_\sigma S\rightarrow
N(S)_\sigma=T_\sigma X/T_\sigma S$,
with $k(Y):=\nabla_YY$(modulo $T_\sigma S$), for $\nabla$ any connection in
the projective structure of $X$. In particular, we take for $X$ the
projective space $\P(T_{\bar\gamma} B)$, with its canonical projective
structure, and for $S$ we take $\P(V^\alpha_{\bar\gamma})$, the
 projectivized $\alpha$-cone in $\bar\gamma$.

\begin{defi}\label{sigma} The projective curvature of the $\alpha$-cone
  $V^\alpha_{\bar\gamma}$ at the generating line $\sigma\subset
  V^\alpha_{\bar\gamma}$ is the projective curvature of
  $\P(V^\alpha_{\bar\gamma})$ in $\sigma$, and is identified to a
  linear application
$$K^\alpha_{\gamma,x}: T_\sigma S\otimes T_\sigma S\rightarrow
N(S)_\sigma,$$
  where $\sigma$ is the tangent direction to $B^\alpha_{\bar\gamma,x}$ in
  $\bar\gamma$, and $S:=\P(V^\alpha_{\bar\gamma})$.
\end{defi}
In order to compute the projective curvature of $V^\alpha_{\bar\gamma}$, we
establish first some canonical isomorphisms between the spaces
appearing in the above definition and some linear subspaces of $T_x\mathbf{M}$.
We will fix now the geodesic $\gamma$, the point $x\in\gamma$
(therefore also $\sigma=T_{\bar\gamma}
B^\alpha_{\bar\gamma,x}\in\P(T_{\bar\gamma} B)$), and, thus, the
$\alpha$-plane $F^\alpha_x\subset T_x\mathbf{M}$ containing $\dot\gamma_x$, as
well as $\dot\gamma^\perp_x\subset T_x \mathbf{M}$, the orthogonal space to 
$\dot\gamma_x$. 

For simplicity, in the following lemmas we will currently omit some
indices referring to these fixed objects.

\begin{lem}\label{l1}
  There is a canonical isomorphism $\tau$ between the tangent space
  $T_\sigma S$ to the projective cone $S=\P(V^\alpha_{\bar\gamma})$ and the
  tangent space $T_x\gamma$ to the geodesic $\gamma$ in the point $x$
  corresponding to the direction $\sigma\in\P(T_{\bar\gamma} B)$. 
\end{lem}
\begin{proof} Let $Y\in T_x\gamma$. We will define $\tau^{-1}(Y)$ as
  follows: Recall that $T_\sigma S\simeq\mbox{Hom}(\sigma,E/\sigma)$,
  where $E(=E_x):=T_\sigma V^\alpha_{\bar\gamma}$ (the tangent space in a
  point to a cone depends only on the line containing the point). We
  know that $\sigma$ corresponds to $\j^\alpha_{\gamma,x}$, the space
  of  Jacobi fields on $\gamma$,
  vanishing at $x$, and such that $\dot J_x\in F^\alpha$. It will be
  shown in the proof of the next theorem that $E$ consists of classes
  of Jacobi vector fields such that $J_x,\dot J_x\in F^\alpha$,
  (\ref{*w}).

Then, on a representative Jacobi field $J\in \j^\alpha_{\gamma,x}$,
we define $\tau^{-1}(Y)$ to be the class of Jacobi fields in
$E/\sigma$, represented by the following Jacobi field $J^Y$ on
$\gamma$, which is given by $J^Y_x:=\nb_Y J, \; \dot J^Y_x:=0$. We
remark that $\nb_Y J$ is what we usually note $\dot J$, when the
parameter on $\gamma$ is understood.

It is straightforward to check that $J\mapsto J^Y$ induces an
isomorphism $\tau^{-1}(Y):\sigma\rightarrow E/\sigma$ for each non-zero
$J\in\sigma=\j^\alpha_{\gamma,x}/\j^\gamma_\gamma$.
\end{proof}

We remark that $V^\alpha_{\bar\gamma}\subset V_{\bar\gamma}$, the
4-dimensional 
subspace represented by Jacobi fields $J$, such that $J,\dot
J\perp\dot\gamma$. We further introduce the subspace
$H^\alpha_{\bar\gamma,x}\subset V_{\bar\gamma}$, represented by Jacobi fields
$J$ as before, with the additional condition $J_x\in F^\alpha_x$. It
is a 3-dimensional subspace, and it contains $E_x$. The curvature of
$V^\alpha_{\bar\gamma}$ will take values in $\mbox{Hom}(TS\otimes TS,N^V(S))$,
and we will show (\ref{**w}) in the proof of the next theorem that it
takes  values
in a smaller space, $\mbox{Hom}(TS\otimes TS,N^H(S))$. $N_\sigma^V(S)$ is
just the normal space of $S$ in $\P(V_\gamma)$ at $\sigma$, and
$N^H(S)$ is the 
subspace of $N^V_\sigma(S)$ consisting in elements represented by
$\xi\in\mbox{Hom}(\sigma,H^\alpha_{\bar\gamma,x})
\subset\mbox{Hom}(\sigma,V_{\bar\gamma})$. 

\begin{lem}\label{l2}
  There is a canonical isomorphism $$\rho:N^H(S)\rightarrow
  \mathrm{Hom}(F^\alpha/T\gamma, \gamma^\perp/F^\alpha).$$
\end{lem}
\begin{proof} As $H$ is a subbundle of the normal bundle $N(S)$, $N^H(S)$ is
  isomorphic to\linebreak $\mbox{Hom}(\sigma,H/E)$. As in Lemma \ref{l1}, we
  will construct the inverse isomorphism $\rho^{-1}$: Let 
  $\xi:F^\alpha/T\gamma\rightarrow \gamma^{\perp}/F^{\alpha}$ be a
  linear application. Let $\xi_0:F^\alpha\rightarrow \gamma^\perp$ be
  a representant of $\xi$ (it involves a choice of a complementary
  space to $F^\alpha$ in $\gamma^\perp$). We define $\rho^{-1}(\xi)\in
  \mbox{Hom}(\sigma,H/E)$ as being induced by the following linear
  application between spaces of Jacobi fields on $\gamma$: 

$\rho^{-1}(\xi):\j^\alpha_{\gamma,x}\rightarrow
\j^{\alpha,\perp}_{\gamma,x}$, where the second space corresponds to
$H_x$, i.e. it contains Jacobi fields $J$ such that $J_x\in F^\alpha,
\dot J_x\perp \dot\gamma_x$. Consider a parameterization of $\gamma$
around $x$, and let $J\in\j^\alpha_{\gamma,x}$. We define
$J^\xi:=\rho^{-1}(\xi) (J)$ by $J^\xi_x:=0,\dot J^\xi_x:=\xi_0(\dot
J_x)$, and it is easy to check that the class of $J^\xi$ in $H/E$ is
independent  of the representant $\xi_0$, such that $\rho^{-1}$ is
well-defined. It is also obviously invertible.
\end{proof}

We are now in position to translate the projective curvature of
$V^\alpha_\gamma$ in terms of conformal invariants of $(\mathbf{M},c)$:

\begin{thrm}\label{t1}
  Let $x$ be  a point in a null-geodesic $\gamma$. Then the projective curvature
  $K$ of the $\alpha$-cone $V^\alpha_{\bar\gamma}$ at $\sigma$
  (corresponding to $x$, see Definition \ref{sigma}), which is a
  linear map 
$$K:T_\sigma S\otimes T_\sigma S\rightarrow N^V(S)_\sigma,$$
takes values in $N^H(S)_\sigma$ (see above), and
is canonically identified to the linear map
$$K':T_x\gamma\otimes T_x\gamma\rightarrow
\mathrm{Hom}(F^\alpha_x/T_x\gamma,\gamma^\perp_x/F^\alpha_x),$$
defined by the self-dual Weyl tensor of {\bf M}:
$$K'(Y,Y)(X)=W^+(Y,X)Y,\; Y\in T_x\gamma,\, X\in F^\alpha_x.$$
\end{thrm}
\begin{proof} 
Consider the following analytic map, which parameterizes, locally
around $x\in\gamma$,  the
deformations of the geodesic $\gamma$ that correspond to points
contained in the integral $\alpha$-cone $B^\alpha_{\bar\gamma}$: 
$$f:U\rightarrow \mathbf{M},\;  f(t,s,u)=\gamma^{t,s}(u),$$
where $U$ is a neighborhood of the origin in $\C^3$, and 
$\gamma^{t,s}$ is a deformation of the null-geodesic $\gamma$, such that
$$\gamma^{t,s}(t)=\gamma(t),\; \dot\gamma^{t,s}(t)\in
F^\alpha_{\gamma(t)},$$
where the parameterization of the geodesic $\gamma$ satisfies
$\gamma(0)=x$, and $F^\alpha_{\gamma(u)}$ is the $\alpha$-plane in
$T_{\gamma(u)}\mathbf{M}$ containing $\dot\gamma(u)$. 

{\bf Convention} We know that $f$ is defined around the origin in
$\C^3$, so there exists a polydisc centered in the origin included in
$U$, therefore all the relations that we will use are true for
values of the variables $t,s,u$ sufficiently close to 0. For
simplicity, we will omit to mention these domains.
\medskip 

The geodesics $\gamma^{t,s}$ correspond to points in
$B^\alpha_{\bar\gamma,\gamma(t)}$, and the Jacobi fields $J^t$ on
$\gamma$, defined as 
$$J^t(u):=\d s f(t,0,u)\in T_{\gamma(u)}\mathbf{M},$$
correspond to vectors in $V^\alpha_{\bar\gamma}$ tangent to the above
mentioned lines. We suppose that the deformation $f$ is {\it
  effective}, i.e. $\d u \gamma^{t,s}(u)\ne 0$ and $J^t\not\in
\j^\gamma_\gamma $, which is equivalent to $\dot J^t(t)\not\in
T_{\gamma(t)}\gamma $. In order to compute the projective curvature of
$V^\alpha_{\bar\gamma}$, we need thus to study the (second order)
infinitesimal variation of these Jacobi fields on $\gamma$. As they
are determined by their value and first order derivative in
$\gamma(0)=x$, we need to evaluate $\d t J^t(0)|_{t=o},\d t\dot
J^t(0)|_{t=0}$  for
the first derivative of $J^t$ at $t=0$, and $\d t^2 J^t(0)|_{t=0},
\d t^2 \nb \dot J^t(0)|_{t=0}$ for the second. Dots mean, as before,
covariant differentiation with respect to the ``speed'' vector
$\dot\gamma$, thus correspond to the operator $\d u$.

As the covariant derivation $\nb$ has no torsion, we can apply the
usual commutativity relations between the operators $\d t, \d s, \d t$
and use them to differentiate the following equation, which follows
directly from the definition of $f$ and $J^t$:
\begin{equation}\label{eq:jt}
J^t(t)=0 \  \forall t.
\end{equation}
We get then
\begin{equation}
  \label{eq:djt}
  \d t J^t(t)+\dot J^t(t)=0, 
\end{equation}
We recall now that, besides (\ref{eq:jt}), we have $\dot J^t(t)\in
F^\alpha_{\gamma(t)}$, thus $\dot J^t(t)$ is isotropic, which implies that:
\begin{equation}
  \label{eq:djt2}
  \langle \d t \dot J^t(t),\dot J^t(t)\rangle =0,
\end{equation}
as $\ddot J^t(t)=R(\dot\gamma(t),J^t(t))\dot\gamma(t)=0$. Equations
(\ref{eq:djt}) and (\ref{eq:djt2}) prove that
\begin{equation}\label{*w} \d t J^t|_{t=0}\in\j^\alpha_{\gamma,x},
\end{equation}
which completes the proof of Lemma \ref{l1}. From (\ref{eq:djt2}), it
equally follows that $\d t
\dot J^t(t)$ is isotropic, and, by differentiating (\ref{eq:djt2}), we
get 
\begin{equation}
  \label{eq:djt3}
  \langle \d t^2 \dot J^t(t),\dot J^t(t)\rangle =-\langle \d t\ddot J^t(t),\dot J^t(t)\rangle .
\end{equation}

From, (\ref{eq:djt}) we have that $\d t J^t(t)$ is isotropic, and also
$$\d t^2 J^t(t)+2\d t\dot J^t(t)=0,$$
which, together with (\ref{eq:djt2}), implies that $\d t^2
J^t(0)|_{t=0} \in F^\alpha_x$. We have then
\begin{equation}
  \label{**w}
  \d t^2 J^t|_{t=0}\in\j^{\alpha,\perp}_{\gamma,x},
\end{equation}
which proves that the curvature $K$ of the $\alpha$-cone takes values
in $N^H(S)$, as it is represented by $\d t^2 J^t|_{t=0}$.

In view of the Lemmas \ref{l1} and \ref{l2}, it is clear now that the
projective curvature $K$ is represented by the following application:
$$(\dot\gamma,\dot\gamma,\dot J)_x\longmapsto \d t^2 J^t(0)|_{t=0}.$$
From (\ref{eq:djt3}), as $\d t\ddot J^t(t)=R(\dot\gamma,\d
tJ^t)\dot\gamma$ and $\d tJ^t(t)=-\dot J^t(t)$, we get
$$\langle K(\dot\gamma,\dot\gamma)(\dot J),\dot J\rangle =\langle R(\dot\gamma,\dot
J)\dot\gamma,\dot J\rangle .$$

The right hand side actually involves only $W^+$, as the other
components of the Riemannian curvature vanish on this combination of
vectors, thus we can replace $R$ with $W^+$ inn the above relation. On
the other  hand, the class of $W^+(\dot\gamma,\dot J)\dot\gamma$
modulo $F^\alpha$ is determined by its scalar product with $\dot J$,
which represents a non-zero generator of $F^\alpha/T\gamma$.

The proof of the Theorem is now complete.
\end{proof}

\obs We may ask whether the projective lines in $Z$ are the geodesics
of some projective structure. Indeed, in the conformally flat case,
when {\bf M} is the Grassmannian of 2-planes in $\C^4$ (the
complexification of the Moebius 4-sphere), $Z\simeq \cp3$, and the
complex lines are given by the standard (flat) projective
structure. But there are two reasons (related to each other, as we
will soon see) why $Z$ cannot carry a canonical projective structure:
First, we do not necessarily have projective lines $Z_x\ni\bar\beta$
in every direction of $T_{\bar\beta} Z$ (this would mean that
$\beta\simeq\cp2$, see next Section for a treatment of this problem),
and second, the lift of a 2-plane $F^\gamma\subset T_{\bar\beta} Z$
would be a 2-plane in $T_{\bar\gamma} B$, so $V^\alpha_{\bar\gamma}$
would be a flat cone: 

\begin{corr}\label{twfl}
  The projective lines $Z_x$ in the twistor space $Z$ are geodesics of
  a projective structure iff it is projectively flat, and {\bf M} is
  conformally flat. 
\end{corr}
\begin{proof} If $Z$ admits a projective structure, some of whose
  geodesics are the lines $Z_x$, then we have, for a fixed $\bar\beta\in
  Z$, a linear connection around $\bar\beta$, whose geodesics in the
  directions of $Z_x,\,\bar\beta\in Z_x\,(\Leftrightarrow x\in\beta\subset
  \mathbf{M})$ coincide, locally, with $Z_x$. This means that the integral
  $\alpha$-cone $Z^\gamma$, for $\gamma\subset\beta$ a null-geodesic, is part
  of a complex surface (namely $\exp(F^\gamma)$, where
  $F^\gamma\subset T_{\bar\beta} Z$ is the 2-plane corresponding to
  $\gamma$). Then the integral $\alpha$-cone $B^\alpha_{\bar\gamma}$, the
  lift to $B$ of $Z^\gamma$, is also a complex surface, thus
  $V^\alpha_{\bar\gamma}$ is a subset of the tangent space
  $T_{\bar\gamma}B^\alpha_{\bar\gamma}$, thus a flat cone. As this is
  true for all points of $Z$ and for all null-geodesics $\gamma$, Theorem
  \ref{t1} implies that {\bf M} is flat.

On the other hand, it is well-known that the twistor space of a
conformally flat manifold admits a flat projective structure, for
which the projective lines $Z_x$ are geodesics, \cite{ahs}.
\end{proof}


\section{Compactness of null-geodesics and conformal flatness}\label{s.w5}
\subsection{Complete $\alpha$-cones in $Z$}\label{ss.w51}
We have given, in the preceding Section, a way to measure the
projective curvature of the $\alpha$-cone in $B$; we shall see now
what happens in the special case when this cone is {\it complete} in
a point $\bar\gamma$, i.e. when $\P(V^\alpha_{\bar\gamma})$ is a compact
submanifold in $\P(T_\gamma B)$.

This situation appears for example if, for any direction in
$F^\gamma\subset T_{\bar\beta} Z$, there are projective lines in $Z$
tangent to it.

\begin{thrm}\label{twist}
 Let $Z$ be the twistor space of the connected civilized self-dual 
4-manifold $(\mathbf{M},c)$, and suppose that, for a point $\beta\in Z$ and for 
a 2-plane $F^\alpha\subset T_{\bar\beta} Z$, there are projective lines $Z_x$ 
tangent to each direction of $F^\alpha$. Then $(\mathbf{M},c)$ is conformally flat.
\end{thrm}

\begin{proof} The idea is to prove that the integral $\alpha$-cone 
$Z^\gamma$ is a smooth surface. We know that this holds in all its 
points except for the vertex $\bar\beta$ (Proposition \ref{p1}). 
The fact that all direction in $F^\gamma$ admits a tangent 
line is a necessary condition for this cone to be a smooth surface, as
it needs to be well-defined around $\bar\beta$.

We choose an auxiliary hermitian (real) metric  $h$ on $Z$. Its
restrictions $h_x$ to the lines 
$Z_x\subset Z^\gamma$ yield K{\"a}hlerian metrics on these lines; in fact
these metrics are deformations of one another, just like the lines $Z_x$
are. This means that the metrics $h_x$ depend continuously on $x\in
\P(F^\alpha)$, a parameter in a compact set. We can therefore find a lower 
bound $r_0>0$ for the injectivity radius of all $(Z_x,h_x)$ at
$\bar\beta$, and a finite upper bound $R$ for the norm of all the
second fundamental forms $H_x: TZ_x\otimes TZ_x\rightarrow
(TZ_x)^\perp\, (\subset TZ).$ We can also suppose that $r_0$ is
smaller than the injectivity radius of $(Z,h)$ at $\bar\beta$. 

The first step is to prove that $Z^\gamma$ is a submanifold of class 
$\mathcal{C}^1$. As its tangent space is everywhere a complex subspace of 
$TZ$, it will follow that it is a complex analytic submanifold.

Consider now the exponential map $\exp_{\bar\beta}:T_{\bar\beta}
Z\rightarrow Z$, defined for the metric $h$; If we restrict it to a
ball of radius less than $r_0$, it is a diffeomorphism into $Z$. The
image of the complex plane $F^\alpha$ is then a smooth 4-dimensional
real submanifold $S$ of $Z$, and  there exists a positive number $r_1$
such that the exponential map in the directions normal to $S$,
$$\exp_S:TS^\perp\rightarrow Z,\, \exp (Y):=\exp_y(Y),\,
\mbox{ for }Y\in T_yS^\perp,$$
restricted to the vectors of length less than $r_1$, is a diffeomorphism.

The image of this diffeomorphism is a tubular neighborhood of $S$, and we 
will denote by $N(S,r)$ such a tubular neighborhood of ``width'' $r$,
for $r<r_1$.

The existence of an upper bound $R$ for the second fundamental forms of 
$Z_x,\forall x\in\gamma$ implies the following fact:

\begin{lem} For any $r<r_1$, there is a neighborhood $U\subset
  T_{\bar\beta} Z$ of the origin such that $\exp (U)\cap Z^\gamma$ is
  contained in $N(S,r)$, and is transverse to the fibers of the
  orthogonal projection $p^S:N(S,r)\rightarrow S,\; p^S(\exp(Y)) :=y$,
  where $Y\in T_yS$. 
\end{lem}
 This is standard if $Z^\gamma$ is a submanifold; but it is also true in our 
case, where $Z^\gamma$ is a union of submanifolds $Z_x$.

Now it is easy to prove that $Z^\gamma$ is a $\mathcal{C}^1$ 
submanifold of $Z$ (the projection $p^S$ yields a local $\mathcal{C}^1$ 
diffeomorphism from a neighborhood of $\bar\beta$ in $S$ to a
neighborhood of $\bar\beta$ in $Z^\gamma$; it is $\mathcal{C}^1$ in
$\bar\beta$ because $S$ is tangent to $Z^\gamma$ at $\bar\beta$). 

So $Z^\gamma$ is a $\mathcal{C}^1$ submanifold of $Z$; Its tangent
space is complex in each point, thus $Z^\gamma$ is a complex-analytic
surface immersed in $Z$.

We have then that $B^\alpha_{\bar\gamma}\subset B=\P(T^*Z)$, being the
lift of $Z^\gamma$, is a smooth analytic surface immersed in $B$, in
particular the $\alpha$-cone $V^\alpha_{\bar\gamma}$ is a complex
plane.

Theorem \ref{t1} implies that $W^+$ vanishes on the $\alpha$-plane
$F_x^\alpha\subset T_x\mathbf{M}$ which contains $\dot\gamma_x$, for
every point $x\in\gamma$. Now, the plane $F^\gamma\subset T_{\bar\beta} Z$
is not the only one admitting projective lines $Z_x$ tangent to any of
its directions:
all planes ``close'' to $F^\gamma$ have the same property. Then $W^+$
vanishes on a neighborhood of $\gamma$, 
hence on the whole connected manifold {\bf M}. 
\end{proof}
\obs There is a more general situation where the integral
$\alpha$-cone $Z^\gamma$ through $\beta$ is smooth in $\beta$:
\begin{thmm}\label{3'}
  Suppose that, for each direction $\sigma\in\P(T_\beta Z)$, there is
  a smooth  (non-necessarily compact) curve $Z_\sigma$ tangent to
  $\sigma$, such that :

(i) if $\sigma$ is tangent to a  projective line $Z_x$, then
$Z_\sigma=Z_x$;

(ii) $Z_\sigma$ varies smoothly with $\sigma\in\P(F^\gamma)$.

Then 
$$\bar Z_\beta^\gamma:=\bigcup_{\sigma\in\P(F^\gamma)}Z_\sigma$$ is a smooth
surface around $\beta$, containing the $\alpha$-cone $Z^\gamma$ and
$W^+(F^\gamma_x)=0,\,\forall x\in\gamma$, where $F^\gamma_x\subset
T_x \mathbf{M}$ is the $\alpha$-plane containing $\dot\gamma$.
\end{thmm}
The proof is similar to the one of the previous theorem. Note that, if
there is a direction $\sigma$ which is not tangent to a projective
line $Z_x$, we cannot apply the deformation argument in Theorem \ref{twist} to
conclude that $W^+$ vanishes everywhere.

{\noindent\bf Example. } If $\mathbf{M}=\M$, then $Z=\mathcal{F}$ and there are
some particular planes for which the conditions in Theorem \ref{3'} are
satisfied, although Theorem \ref{twist} never applies to $Z$: for a generic
2-plane $F^\gamma$, the $\alpha$-cone $V^\alpha_\gamma$ is not
flat. The above mentioned particular planes in $TZ$ correspond to the
vanishing of $W^+$ on some particular $\alpha$-planes, but {\bf M} is not
anti-self-dual (see Sections \ref{ss.w83}, and also \ref{ss.w87}, \ref{ss.w88}).

\subsection{Compact, simply-connected null-geodesics in {\bf M}}\label{ss.w52}
Theorem \ref{twist} suggests that the existence of a compact null-geodesic
diffeomorphic to $\cp1$ yields strong constraints upon the conformal
structure of {\bf M}. In fact, we have:

 \begin{thrm}\label{compact}
   If a connected complex self-dual 4-manifold $(\mathbf{M},c)$ admits a compact null-geodesic
   diffeomorphic to $\cp1$, the conformal structure of {\bf M} is flat.
 \end{thrm}
\obs A similar result has been proven by Y.-G. Ye using algebraic
geometry techniques~: if a {\it projective} complex manifoldadmits a
conformal structure having a compact null-geodesic diffeomorphic to
$\cp1$, it is conformally flat \cite{ye}. Note that we do not need {\bf M} to
be compact in Theorem \ref{compact}~; on the other hand, we assume it
to be self-dual.
 \begin{proof} 
We first remark that the main difficulty is the definition of $B$, the
space of null-geodesics , and of $Z$, the twistor space of $(\mathbf{M},c)$, as {\bf M} is
not necessarily civilized. This is only possible on small open sets,
but, in general,  we can not expect to have  any global construction
of this kind. Thus, things that were almost obvious in the twistorial
framework (like the existence of compact deformations of the null-geodesic
$\gamma$), seem much more difficult to prove directly. The idea is to
prove that all null-geodesics close to $\gamma$ are diffeomorphic to
$\cp1$. Then, we show that, conversely, every projective line which is
a deformation of $\gamma$ as a compact curve is a null-geodesic. In particular,
sections in the normal bundle $N(\gamma)$ are induced by (local)
Jacobi fields. We obtain then directly that $W^+=0$.

   \begin{prop}\label{p5}
     Let $\gamma$ be an immersed null-geodesic, diffeomorphic to a projective line
     $\cp1$. Then any local Jacobi field $J$ with $\dot
     J\perp\dot\gamma$ induces a global normal field $\nu^J$ on
     $\gamma$.
   \end{prop}
\begin{proof} If $\gamma$ is just a compact geodesic, it may have
     points of self-intersection, but it is always an immersed
     curve. It is more convenient then to think of $\gamma$ as a
     projective line immersed in {\bf M} rather than the image of this
     immersion. The tangent, normal bundles, etc. are also to be
     thought as bundles over this projective line, still denoted by
     $\gamma$. Tubular neighborhoods of $\gamma$ are then
     neighborhoods of the zero section in the normal bundle
     $N(\gamma)$, small enough to be immersed (non-injectively)  in
     {\bf M} as a neighborhood of the image of $\gamma$.

We first notice that $\gamma$ may be decomposed in the
     union of a two open sets $U_1\cup U_2$, both biholomorphic to the
     unit disk in $\C$, and such that $U_1\cap U_2$ is connected.

Then, for any local metric in $c$, we have a Jacobi equation around a
point $x\in\gamma$, and a Jacobi field $J$ corresponding to prescribed
$J_x,\dot J_x$. It is easy to prove that the (local) normal field
induced by $J$ is independent of the chosen metric. Moreover, this
normal field is the unique solution, for the prescribed 1-jet in $x$
induced by $J_x,\dot J_x$, of a second order differential equation on
$N(\gamma)$:

\begin{lem}
  The Jacobi equations for null-geodesic induce a second order linear
  differential operator $P$ on $N(\gamma)$, depending only on the
  conformal structure $c$ of {\bf M}.
\end{lem}
\begin{proof} For a Levi-Civita connection $\nb$ of a local metric on
  {\bf M}, we locally define the following differential operator on
  $T\mathbf{M}|_\gamma$:
$$P:\Gamma(T\mathbf{M}|_\gamma\otimes
S^2(T\gamma))\rightarrow\Gamma(T\mathbf{M}|_\gamma),$$
by $P(Y,X,X):=\nb_X\nb_XY-\nb_{\nb_XX}Y-R(X,Y)X$. It obviously induces
a (local) differential operator on $N(\gamma)$, and all we need to
show is that, for a different connection $\nb'$, the corresponding
operator $P'$ induces the same one on $N(\gamma)$. First we write
$$P(Y,X,X)=\nb_X[X,Y]+\nb_{[X,Y]}X-[\nb_XX,Y],$$ 
then we recall that another Levi-Civita connection $\nb'$ is related
to $\nb$ by the formula \cite{gconf}:
$$\nb'_AB=\nb_AB+\theta(A)B+\theta(B)A-g(A,B)\theta^\sharp, \mbox{ for
  } \theta\in\Lambda^1\mathbf{M},$$
so we directly obtain:
$$P'(Y,X,X)-P(Y,X,X)=2(\nb_Y\theta)(X)X+2\theta(\nb_XY),$$
thus they induce the same operator on $N(\gamma)$. This one is,
therefore, globally defined (the topology of $\gamma$ is not
important).
\end{proof}

Now, for any $x\in\gamma$, we have an unique solution $\nu$ of $P$, for a
prescribed {\it 1-jet} $j^1(\nu)_x$ (which consists in the values in
$x$ of $\nu$ and of his first-order derivative), globally defined on
every contractible open set $U\ni x$ in $\gamma$. This is because on
any such contractible set the equation $P\nu=0$ becomes a second order
ordinary linear equation on a disk in $\C$, which admits global
holomorphic solutions (unique if we fix the initial conditions).

Take now $x\in U_1\cap U_2$. Then, the two solutions $\nu_1$ and
$\nu_2$, defined on $U_1$, resp. $U_2$, coincide on the connected $
U_1\cap U_2$, so they yield a global solution $\nu$ with the prescribed
initial conditions in $x$.

In particular, this solution is a global section of the normal bundle
$N(\gamma)$.
\end{proof}

After the infinitesimal result, the local one:

\begin{prop}\label{tub}
  Small deformations of $\gamma$ are also compact immersed projective
  lines.
\end{prop}
\obs The tubular neighborhoods considered below are always seen as
images, by a local diffeomorphism, of subsets --- which are,
generally, fiber bundles over $\cp1$ --- of the normal bundle of
$\gamma$, resp. $\tilde\gamma$, the lift of $\gamma$ to $\P(C)$. We
need this because of the possible self-intersections of
$\gamma$; $\tilde\gamma$ is always embedded.

\begin{proof} Consider an auxiliary hermitian metric $h$ on {\bf M}.Then
  $h|_\gamma$ induces the same topology like a round metric $h_0$ on
  the sphere $S^2$. We can define a tubular neighborhood $N(r_1)$ of
  $\gamma$, as the open set of points $y\in \mathbf{M}$ with
  $d(y,\gamma)<r_1$. We choose $r_1$ small enough for $N(r_1)$ to be a
  fiber bundle over $\gamma$ (the fiber $N(r_1)_x$, for $x\in\gamma$,
  being the image of the real 4-plane of $T_x\mathbf{M}$, $h$-orthogonal to
  $T_x\gamma$, by the exponential of $h$).

Take now a finite number of contractible open sets whose union
  covers $\gamma$, and we choose holomorphic metrics in $c$ on each of
  this open sets. Then, on these sets we have connections,
  and it is well-known that any point in such a set has a basis
  of geodesically connected neighborhoods \cite{kn}, \cite{ehr}. We
  choose a  finite number
  of such geodesically convex open sets, that cover $\gamma$, and such
  that they are all included in $N(r_1)$. We have then
$$N(r_1)\subset\bigcup_{i=1}^n U_i\supset \gamma.$$ (Of course, they are
geodesically convex only with respect to some particular metric in
$c$, but we are interested only in the implications involving the null-geodesics,
which are independent of the metric.)

It is immediate, \cite{leb5}, that a geodesically convex set $U_i$ has the
following property: all maximal geodesics are closed submanifolds of
$U_i$ and are contractible as topological spaces. This is particularly
true for null-geodesics included in $U_i$. 

We refine now the covering by another one, $U'_i\subset \bar U'_i\subset U_i$,
such that $U'_i=N(r_2)|_{V_i}$, where $V_i:=U_i'\cap\gamma$, such that
$\gamma\subset \cup_{i=1}^n U'_i$. In fact, 
we ask for $U'_i$ to be restrictions to $V_i$ of a tubular neighborhood 
$N(r_2)$, for $r_2$ sufficiently small ($U'_i$ are ``cylindrical''
neighborhoods). We can easily imagine how to find such a refinement of
the initial covering. 

The proof of the proposition now follows two ideas: first, we consider
a very  
special covering by disks (for the round metric $h_0$) of $\gamma$;
second, we define a neighborhood of $\tilde\gamma$ in $\P(C)$ of
isotropic directions for which, using the special covering of
$\gamma$, we can extend the associated null-geodesics and eventually get compact
ones. We call $\tilde\gamma$ the canonical lift of a null-geodesic in $\P(C)$;
it is an integral curve of the geodesic distribution of lines in
$\P(C)$. 

The first idea has nothing to do with complex analysis; it is just a matter of
metric topology on the round sphere $(S^2,\mbox{can})$.

\begin{lem}\label{ldisc} Let $\{V_i\}_{i=\overline{1,n}}$ be a
  covering of  $(S^2,\mbox{can})\simeq (\gamma,h_0)$  by
  open sets. Then there exists a positive number $r_0$ and a finite
  set of points $x_j\in \gamma$, for $j=\overline{1,N}$, such that:

(i) All disks $D(x,10r_0)$ are contained in at least one of the open
sets $V_i$; 

(ii) The disks $D(x_j,r_0)_{j=\overline{1,N}}$ cover $\gamma\simeq S^2$;

(iii) $100r_0<l$, where $l$ is the diameter of $\gamma\simeq S^2$. 
\end{lem}

An important property of this covering is that all sets, as well as
the intersections of a finite number of them, are convex for the round
metric, thus contractible.

We intend to extend null-geodesics which can be projected diffeomorphically onto
$\gamma$ by means of the fiber projection $\rho:N(r_1)\rightarrow
\gamma$. We will do that step-by-step, extending it over disks $D(x,R)$ 
of increasing radius. But before that, we need to restrict ourselves to some
particularly ``close to $\gamma$'' null directions. 

We need two things: the extensions need to remain within $N(r_2)$, and 
they should also be transverse to the fibers of $\rho$, otherwise the
projection into $\gamma$ would not be an invertible diffeomorphism.

First, we consider the following compact subset of $\P(T\mathbf{M})$:
$$S:=\{L\subset T_yM| y\in \overline {N(r_1)},\, 
L\mbox{ tangent to the fibers of } \rho\},$$
where we say that a complex line $L$ is tangent to a real manifold 
$\rho^{-1}(x),\, x\in\gamma$ if it contains a non-zero (real) vector tangent
to this real submanifold.

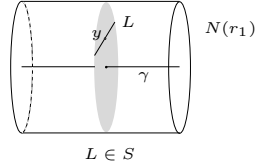
\begin{center}\begin{picture}(0,0)%
\epsfig{file=weyl/tub.pstex}%
\end{picture}%
\setlength{\unitlength}{0.00013100in}%
\begingroup\makeatletter\ifx\SetFigFont\undefined
\def\x#1#2#3#4#5#6#7\relax{\def\x{#1#2#3#4#5#6}}%
\expandafter\x\fmtname xxxxxx\relax \def\y{splain}%
\ifx\x\y   
\gdef\SetFigFont#1#2#3{%
  \ifnum #1<17\tiny\else \ifnum #1<20\small\else
  \ifnum #1<24\normalsize\else \ifnum #1<29\large\else
  \ifnum #1<34\Large\else \ifnum #1<41\LARGE\else
     \huge\fi\fi\fi\fi\fi\fi
  \csname #3\endcsname}%
\else
\gdef\SetFigFont#1#2#3{\begingroup
  \count@#1\relax \ifnum 25<\count@\count@25\fi
  \def\x{\endgroup\@setsize\SetFigFont{#2pt}}%
  \expandafter\x
    \csname \romannumeral\the\count@ pt\expandafter\endcsname
    \csname @\romannumeral\the\count@ pt\endcsname
  \csname #3\endcsname}%
\fi
\fi\endgroup
\begin{picture}(7448,6327)(2483,-7141)
\put(6766,-2011){\makebox(0,0)[lb]{\smash{\SetFigFont{6}{7.2}{rm}
$L$
}}}
\put(5626,-2386){\makebox(0,0)[lb]{\smash{\SetFigFont{6}{7.2}{rm}
$y$
}}}
\put(7411,-4021){\makebox(0,0)[lb]{\smash{\SetFigFont{6}{7.2}{rm}
$\gamma$
}}}
\put(9931,-2251){\makebox(0,0)[lb]{\smash{\SetFigFont{6}{7.2}{rm}
$N(r_1)$
}}}
\put(5356,-7051){\makebox(0,0)[lb]{\smash{\SetFigFont{6}{7.2}{rm}
$L\in S$
}}}
\end{picture}
\end{center}

The hermitian metric $h$ on {\bf M} induces a metric on $\P(T\mathbf{M})$, and also
one on $\P(C)$. We can, then, evaluate the distance between $\tilde\gamma$
and $S$: 
$$\mu_0:=d(\tilde\gamma,S)>0, $$
as they are disjoint compact sets.

Following LeBrun \cite{leb5}, we can define the complex 5-manifolds
$B^i$ as the spaces of null-geodesics of $(U_i,c)$, equivalently the space of
integral curves of the geodesic distribution in $\P(C)|_{U_i}$. The
projections $p_i: \P(C)|_{U_i}\rightarrow B^i$, which send an
isotropic direction to the null-geodesic tangent to it, is a submersion and the
(closed) fibers are precisely the lifts of the null-geodesics of $U_i$. This
construction is possible because $U_i$ are geodesically 
convex (for a particular local holomorphic metric in $c$), see
\cite{leb5} for  details.

We first consider a tubular neighborhood $\tilde N(r^0)$ of
$\tilde\gamma$ in $\P(C)$ which projects, by $\pi:\P(C)\rightarrow \mathbf{M}$
{\it inside} $N(r_2)$, and such that $100r^0<\mu_0$. This second
condition ensures that all directions in $\tilde N(r^0)$ are transverse
to the fibers of $\rho$.  

Consider then the following neighborhoods of $\tilde\gamma|_{\bar
  V_i}$: $p_i:\P(C)|_{U_i}\rightarrow B^i$ is an open application, so
  we define $C_i$ to be $p_i^{-1}(C_B^i)$, where $C_B^i$ is an open
  neighborhood of $\bar\gamma\in B^i$ contained in $p_i(\tilde
  N(r_0))$. 

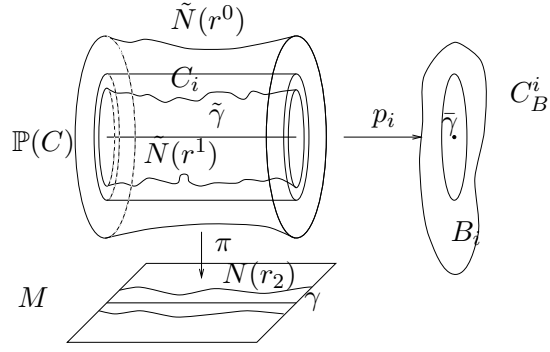
\begin{center}\input{weyl/tubes.pstex_t}\end{center}

It is important to note that $C_i$ have the following
  property:

  \begin{lem}
    for any point $Y\in C_i$, the null-geodesic $\gamma^Y$ tangent to $Y$,
    contained in
    $U_i$, lies into $N(r_2)$ and is always transverse to the fibers
    of $\rho:N(r_1)\rightarrow \gamma$. Hence, its restriction to the
    points of $U'_i$ projects diffeomorphically onto
    $V_i\subset\gamma$. Moreover, all the points of $\tilde\gamma^Y$
    that lie over $\bar V_i$ are in $C_i$.
  \end{lem}

Obviously, the crucial property of these open sets is that every null-geodesic
starting there is totally contained in $C_i$, at least the part that
``lies over'' (in the sense of the projection $\rho$) $\bar
V_i$. After constructing these sets $C_i,\forall i$, we define $r^1>0$
small enough for the tubular neighborhood $\tilde N(r^1)$, restricted
to $\bar V_i$, to be contained into $C_i$, for all $i$. For each $i$,
this means that $r^1$ has to be less than the minimum of the following
continuous functions defined on the compact set $\bar V_i$:
$$\bar V_i\ni x\mapsto d(T_x\gamma, \P(C)_x\smallsetminus C_i).$$

This neighborhood $\tilde N(r_1)$ of $\tilde\gamma$ has the following
property:
\begin{lem}
  For each $Y\in \tilde N(r_1)$, and for each $i$ such that
  $\pi(Y)=x\in V_i$, we have $Y\in C_i$, and thus the whole null-geodesic
  $\gamma^{Y,i}$, contained in $U_i$, is lifted to
  $\tilde\gamma^{Y,i}\subset C_i$.
\end{lem}

The disadvantage of $\tilde N(r_1)$ is that it does not necessarily
contain $\tilde\gamma^{Y,i}$. But we know that the latter is contained in
$\tilde N(r_0)$, which contains the union of all $C_i$. 

We recall now that the idea of proof is to extend a null-geodesic $\gamma^Y$
close  to $\gamma$ over the disks $D_j:=D(x_j,r_0)$. Every extension
over a disk brings $\gamma^Y$ from $\tilde N(r^1)$ to the larger set
$\tilde N(r^0)$. As we have a finite, well-determined, number of disks
$N$, all we need now to apply our extending idea is a sequence of open
sets $\tilde N(r^k)$ such that 
\begin{equation}
  \label{eq:ppp}
  \forall Y\in \tilde N(r^k), \mbox{ such that } x:=\pi(Y)\in V_i,\;
  \tilde\gamma^{Y,i}\subset\tilde N(r^{k-1}).
\end{equation}

To do that, we construct $C^1_i\subset \tilde N(r^1)$ as we have done
for $C_i$, and then $\tilde N(r^2)$ by repeating the same
procedure. We stop after $N$ (the number of disks covering
$(\gamma,h_0)$, see Lemma \ref{ldisc}) steps and claim: 

\begin{prop}
  $\forall Y\in \tilde N(r^N)$, the null-geodesic $\gamma^Y$ extends to a
  compact curve which projects (via $\rho$) diffeomorphically onto
  $\gamma$.
\end{prop}

\begin{proof} Fix $x:=\pi(Y)$. We can define $\gamma^Y$ inside $U_i$,
where $U_i$ is an open set containing $x$.
In particular, $\gamma^Y$ is well-defined over $D(x_0,10r_0)$, where 
$x_0:=\rho (x)$, see Lemma \ref{ldisc}. (Of course, this is because
$\gamma^Y|_{U_i}$ is transverse to the fibers of $\rho$. We intend to
extend it over disks centered in $x_0$. 

Consider the domains $D^x_i:=D^{x_0}_i$ which are ``quadrilaterals''
contained in $D(x_i,10r_i)$ and containing $D(x_i,r_0)$, as in the
following picture (the ``vertical'' parts of the border of $D^x_i$ are
segments of circles centered in $x_0$): 

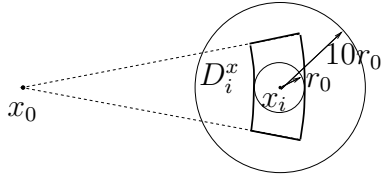
\begin{center}\begin{picture}(0,0)%
\epsfig{file=weyl/dom.pstex}%
\end{picture}%
\setlength{\unitlength}{0.00030600in}%
\begingroup\makeatletter\ifx\SetFigFont\undefined
\def\x#1#2#3#4#5#6#7\relax{\def\x{#1#2#3#4#5#6}}%
\expandafter\x\fmtname xxxxxx\relax \def\y{splain}%
\ifx\x\y   
\gdef\SetFigFont#1#2#3{%
  \ifnum #1<17\tiny\else \ifnum #1<20\small\else
  \ifnum #1<24\normalsize\else \ifnum #1<29\large\else
  \ifnum #1<34\Large\else \ifnum #1<41\LARGE\else
     \huge\fi\fi\fi\fi\fi\fi
  \csname #3\endcsname}%
\else
\gdef\SetFigFont#1#2#3{\begingroup
  \count@#1\relax \ifnum 25<\count@\count@25\fi
  \def\x{\endgroup\@setsize\SetFigFont{#2pt}}%
  \expandafter\x
    \csname \romannumeral\the\count@ pt\expandafter\endcsname
    \csname @\romannumeral\the\count@ pt\endcsname
  \csname #3\endcsname}%
\fi
\fi\endgroup
\begin{picture}(6083,2926)(901,-3114)
\put(6016,-1591){\makebox(0,0)[lb]{\smash{\SetFigFont{12}{14.4}{rm}
$r_0$
}}}
\put(6346,-1201){\makebox(0,0)[lb]{\smash{\SetFigFont{12}{14.4}{rm}
$10r_0$
}}}
\put(5251,-1936){\makebox(0,0)[lb]{\smash{\SetFigFont{12}{14.4}{rm}
$x_i$
}}}
\put(901,-2101){\makebox(0,0)[lb]{\smash{\SetFigFont{12}{14.4}{rm}
$x_0$
}}}
\put(4186,-1561){\makebox(0,0)[lb]{\smash{\SetFigFont{12}{14.4}{rm}
$D_i^x$
}}}
\end{picture}
\end{center}

 We change, if necessary, the
order of the indices $i$ of $x_i$ and $D^x_i$ such that it coincides with the 
ordering of increasing distances $d(x_0,x_i)$. We define than the open sets
$$\Delta_k:=\left\{\begin{array}{ll}D(x_0,r_0) & \mbox{ if
      $d(x_0,x_k)<9r_0$,}\\
\Delta_{k-1} & \mbox{ if $d(x_0,x_k)>l-10r_0$,}\\
\Delta_{k-1}\cap D^x_k & \mbox{ otherwise}.\end{array}\right.$$
\obs The closed disk $\bar D(x_0,d(x_0,x_{k+1}-r_0)$ is included in
      $\Delta_k$  as soon as $x_k\not\in D(\bar x_0,10r_0)$, where
      $\bar x_0$ is the point of $\gamma=S^2$ opposed to $x_0$. Then,
      because of the specific geometry of the domains $D^x_i$ (see
      Lemma \ref{ldisc}), we easily conclude that the domains
      $\Delta_k$ are contractible (along the geodesics of the sphere
      passing through $x_0$) and so are the intersections
      $\Delta_k\cap D^x_{k+1}$, too. 

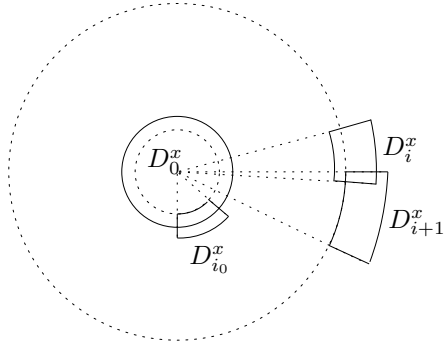
\begin{center}\begin{picture}(0,0)%
\epsfig{file=weyl/disk.pstex}%
\end{picture}%
\setlength{\unitlength}{0.00070000in}%
\begingroup\makeatletter\ifx\SetFigFont\undefined
\def\x#1#2#3#4#5#6#7\relax{\def\x{#1#2#3#4#5#6}}%
\expandafter\x\fmtname xxxxxx\relax \def\y{splain}%
\ifx\x\y   
\gdef\SetFigFont#1#2#3{%
  \ifnum #1<17\tiny\else \ifnum #1<20\small\else
  \ifnum #1<24\normalsize\else \ifnum #1<29\large\else
  \ifnum #1<34\Large\else \ifnum #1<41\LARGE\else
     \huge\fi\fi\fi\fi\fi\fi
  \csname #3\endcsname}%
\else
\gdef\SetFigFont#1#2#3{\begingroup
  \count@#1\relax \ifnum 25<\count@\count@25\fi
  \def\x{\endgroup\@setsize\SetFigFont{#2pt}}%
  \expandafter\x
    \csname \romannumeral\the\count@ pt\expandafter\endcsname
    \csname @\romannumeral\the\count@ pt\endcsname
  \csname #3\endcsname}%
\fi
\fi\endgroup
\begin{picture}(2855,2534)(758,-2003)
\put(2116,-1411){\makebox(0,0)[lb]{\smash{\SetFigFont{10}{12.0}{rm}$D^x_{i_0}$}}}
\put(1801,-691){\makebox(0,0)[lb]{\smash{\SetFigFont{10}{12.0}{rm}$D^x_0$}}}
\put(3556,-601){\makebox(0,0)[lb]{\smash{\SetFigFont{10}{12.0}{rm}$D^x_i$}}}
\put(3601,-1141){\makebox(0,0)[lb]{\smash{\SetFigFont{10}{12.0}{rm}$D^x_{i+1}$}}}
\end{picture}
\end{center}

We prove, then, by induction, that $\gamma^Y$ can be extended over 
$\Delta_k$, and that all the corresponding points of $\tilde\gamma^Y$
are contained in $\tilde N(r^{N-k})$. 

This is obvious for small values of $k$. When we add $D^x_k$ to 
$\Delta_{k-1}$, we consider a point $z$ in the connected (see above)
intersection $ \Delta_{k-1}\cap D^x_k$. It is contained in the $U'_i$
that contains $D(x_k,10r_0)$ and, as $\tilde\gamma^Y|_{\Delta_{k-1}}$
is contained in $\tilde N(r^{N-k+1})$, the connected piece
$\tilde\gamma^Y|_{\Delta_{k-1}\cap D^x_k}$ is contained in $C_i^{N-k}$,
and thus the whole extension $\tilde\gamma^{Y,i}$ in $U_i$ is
contained in $C_i^{N-k}$, hence in $\tilde N(r^{N-k})$. The
connectedness of the considered piece implies that
$\gamma^{Y,i}|_{\Delta_{k-1}\cap D^x_k}$ coincides with
$\tilde\gamma^Y|_{\Delta_{k-1}\cap D^x_k}$. Thus we have obtained an
extension of $\gamma^Y$ to $\Delta_k$, such that the corresponding
lift lies in $\tilde N(r^{N-k})$, as claimed. 

For large values of $k$, the $\Delta_k$ are all identical to 
$\Delta_{N-1}$, which contains $D(x_0,l-10r_0)$. We have thus proven 
that there is an extension of $\gamma^Y$ over this disk, such that
the corresponding points of $\tilde\gamma^Y$ are in $\tilde N(r^1)$.
Consider then the disk $D(\bar x_0,10r_0)$; it is contained in some
$V_i$:
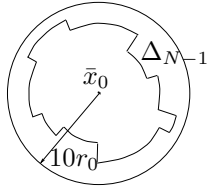
\begin{center}\begin{picture}(0,0)%
\epsfig{file=weyl/deltan.pstex}%
\end{picture}%
\setlength{\unitlength}{0.00017500in}%
\begingroup\makeatletter\ifx\SetFigFont\undefined
\def\x#1#2#3#4#5#6#7\relax{\def\x{#1#2#3#4#5#6}}%
\expandafter\x\fmtname xxxxxx\relax \def\y{splain}%
\ifx\x\y   
\gdef\SetFigFont#1#2#3{%
  \ifnum #1<17\tiny\else \ifnum #1<20\small\else
  \ifnum #1<24\normalsize\else \ifnum #1<29\large\else
  \ifnum #1<34\Large\else \ifnum #1<41\LARGE\else
     \huge\fi\fi\fi\fi\fi\fi
  \csname #3\endcsname}%
\else
\gdef\SetFigFont#1#2#3{\begingroup
  \count@#1\relax \ifnum 25<\count@\count@25\fi
  \def\x{\endgroup\@setsize\SetFigFont{#2pt}}%
  \expandafter\x
    \csname \romannumeral\the\count@ pt\expandafter\endcsname
    \csname @\romannumeral\the\count@ pt\endcsname
  \csname #3\endcsname}%
\fi
\fi\endgroup
\begin{picture}(5484,5484)(4232,-6447)
\put(8429,-2535){\makebox(0,0)[lb]{\smash{\SetFigFont{10}{12.0}{rm}
$\Delta_{N-1}$
}}}
\put(5579,-5850){\makebox(0,0)[lb]{\smash{\SetFigFont{10}{12.0}{rm}
$10r_0$
}}}
\put(6629,-3270){\makebox(0,0)[lb]{\smash{\SetFigFont{10}{12.0}{rm}
$\bar x_0$
}}}
\end{picture}
\end{center}
The intersection $\Delta^x$ of this disk with $\Delta_{N-1}$ 
is a connected open subset of $V_i$, and we know that $\tilde\gamma^Y|_
{\Delta^x}$ is contained in $\tilde N(r^1)$. This implies that we can 
extend in a unique way  $\tilde\gamma^Y|_{\Delta^x}$ to $V_i$, in 
particular to  $D(\bar x_0,10r_0)$, and the corresponding points in 
$\tilde\gamma^Y|_{V_i}$ are in $C_i\subset\tilde N(r^0)$.

We have proven that $\gamma^Y$ extends over $\gamma$, i.e. there is a 
maximal extension (obviously unique) of the null-geodesic tangent to $Y$, such that
it projects (via $\rho$) diffeomorphically onto $\gamma$. The projection
is $C^\infty$, but the extended null-geodesic is clearly an analytic submanifold
of {\bf M}.
\end{proof}

The proof of Proposition \ref{tub} is now complete.
\end{proof}

\obs Generic deformations of $\gamma$
are embedded projective lines. Indeed, let $\{x_1,...,x_k\}$ be the
nodes (self-intersection points) of $\gamma$, and consider the
manifold $\widetilde{\mathbf{M}}$, obtained by blowing-up the points
$\{x_1,...,x_k\}$. Then, generically, any null-geodesic $\gamma'$ close to
$\gamma$ avoids these points, hence is diffeomorphic to its lift to
$\widetilde{\mathbf{M}}$, which is a deformation of the lift of $\gamma$, thus it is
an immersed projective line. But the lift of $\gamma$ is embedded, and
so must be its deformations, hence $\gamma'$ is embedded. 
\smallskip

The next step is to prove that all deformations of $\gamma$ as a
compact curve are null-geodesics , by a dimension-counting argument; we need to
compute the normal bundle of $\gamma$.
\smallskip

We ask now if the family of projective lines in {\bf M} defined as the 
deformations of $\gamma$ is {\it locally complete} in the sense of 
Kodaira \cite{kod}. For this, we need to prove that the dimension
of the space of global sections in $N(\gamma)$ is equal to 5, i.e. 
to the dimension of $B$. The extensions of the null-geodesics close to
$\gamma$ yield local diffeomorphisms between neighborhoods of 
$\gamma$ in $B^i$, resp. $B^j$. In fact, we have a projection
$p:\tilde N(r^N)\rightarrow W^B$ onto the space of integral curves of
the geodesic distribution in $\tilde N(r^N)\subset \P(C)$. $W^B$ is
the space of complex null-geodesic {\it close to $\gamma$} in {\bf M}. But
essential for us is that $p$ is a submersion, fact that has important
consequences for the normal bundle of $\gamma$ in {\bf M}.

\begin{prop}
The normal bundle of $\gamma$ in {\bf M} is isomorphic to 
$\o1\oplus\o1\oplus\o0$. 
\end{prop}
\begin{proof} It is well-known that all holomorphic bundles over
  $\cp1$ are direct sums of line bundles, all of which are isomorphic
  to $\o k,\, k\in\Z$.

We have the subbundle $N^\perp(\gamma)$ of $N(\gamma)$, 
  represented by vectors orthogonal to $\dot\gamma$. We have the
  following exact sequence:
  \begin{equation}
    \label{eq:o0}
    0\rightarrow N^\perp(\gamma)\rightarrow N(\gamma)\rightarrow
    N(\gamma)/N^\perp(\gamma)\rightarrow 0.
  \end{equation}
The right hand term of this sequence is a line bundle, and it admits
global non-zero sections (extensions of Jacobi fields $J$ such that
$\langle J,\dot\gamma\rangle $ is a non-zero constant, see Proposition
\ref{p5}). $N(\gamma)/N^\perp(\gamma)$ is then isomorphic to $\o a$,
with $a\in\N$.

We denote by $N^\beta(\gamma)$ the subbundle of the normal bundle
represented by vectors in $T\beta$. It admits global sections, namely
the extensions of Jacobi fields contained in $T_\beta$, see
Proposition \ref{p5}. It is a line
bundle, thus isomorphic to $\o {c_1},\; c_1\geq 1$, as it contains
global sections with prescribed 1-jet in a point.

\obs In general, $\o k$ is the line bundle over $\cp1$ admitting
global sections for any prescribed $k$-jet in a point $x$. This
section is unique, and it gives a unique value of the $k+1$-jet in
$x$. 
\smallskip

We have the following exact sequence:
\begin{equation}
  \label{eq:o1}
  0\rightarrow N^\beta(\gamma)\rightarrow N^\perp(\gamma)\rightarrow
  N^\perp(\gamma)/N^\beta(\gamma) \rightarrow 0.
\end{equation}
It is easy to check that the right hand term admits local sections
represented by Jacobi fields for any prescribed 1-jet in a point
$x$. Hence, $N^\perp(\gamma)/N^\beta(\gamma)\simeq\o {c_2},\; c_2\geq
1$.

All we can obtain now is that $N^\perp(\gamma)\simeq
\o{b_1}\oplus\o{b_2}$, with $b_1+b_2=c_1+c_2\geq 2$. Actually
$b_1,b_2\geq 1$, otherwise all sections of $N^\perp$ that vanish
somewhere would be contained in a line subbundle, which would
contradict Proposition \ref{p5}.

We have then $N(\gamma)\simeq\o a\oplus\o {b_1}\oplus\o {b_2},\; a\geq
0,\; b_1,b_2\geq 1$. We want to prove that there is equality in all
these three inequalities.

We know, from Proposition \ref{tub},  that there is a tubular
neighborhood $\tilde N(r^N)$ of $\tilde\gamma$ in $\P(C)$ such that the
null-geodesic distribution yields a foliation with compact leaves, and such that
the projection onto the space $W^B$ of these compact curves is a
submersion. (Of course, this space is nothing but the space of complex
null-geodesics {\it close to} $\gamma$.)

It is obvious then that the normal bundle of $\tilde\gamma$ in $\P(C)$
is trivial, as $\tilde\gamma$ is e fiber of a submersion.

We have now the following exact sequence of bundles, related to the
projection $\pi :\P(C)\rightarrow \mathbf{M}$:
\begin{equation}
  \label{eq:opi}
  0\rightarrow N^\pi(\tilde\gamma)\rightarrow
  N(\tilde\gamma)\rightarrow\pi^* N(\gamma)\rightarrow 0,
\end{equation}
where $N^\pi(\tilde\gamma)$ is the normal subbundle of $\tilde\gamma$
represented by vectors tangent to the fibers of $\pi$. In a point
$T_x\gamma\in\tilde\gamma\subset\P(C)$, the fiber of $\pi$ is equal to
$\P(C)_x$, so the tangent space to it is isomorphic to
$\mbox{Hom}(T_x\gamma,T_x\gamma^\perp/T_x\gamma)$, for the projective
variety $\P(C)_x\subset\P(T_x\mathbf{M})$. Thus 
$$N^\pi(\tilde\gamma)\simeq\mbox{Hom}(T\gamma,N^\perp(\gamma))\simeq
\o {-2}\otimes(\o {b_1}\oplus\o{b_2}).$$

The central bundle in the exact sequence (\ref{eq:opi}) is
trivial. The equation above then implies that the Chern number of
$N(\gamma)$ is subject to the following constraint: 
$$a+b_1+b_2+(b_1-2+b_2-2)=0,$$
thus, as $b_1,b_2\geq 1$ and $a\geq 0$, we have $a=0$ and $b_1=b_2=1$.

\end{proof}
As observed above, generic, compact, simply-connected null-geodesics
are embedded. From now on, we suppose $\gamma$ is one of them.

Then $H^1(N(\gamma))=0$, so we can apply the theory of Kodaira to
deform $\gamma$, \cite{kod}, so the dimension of the space of global
sections of $N(\gamma)$ is $1+2+2=5$, the same as the space of complex
null-geodesic close to $\gamma$, hence 

\begin{corr}
  The deformations of $\gamma$ as a compact curve are null-geodesics in $(\mathbf{M},c)$.
\end{corr}

This means that any global section in $N(\gamma)$ can be represented,
locally, by the Jacobi fields that yield the same element in
$J^1N(\gamma)$, the space of jets of order 1 in $N(\gamma)$.

Recall now the exact sequence (\ref{eq:o1}); we conclude that
$$N^\beta(\gamma)\simeq N^\perp(\gamma)/N^\beta(\gamma)\simeq\o1.$$
We have the canonical isomorphism $N^\alpha(\gamma)\rightarrow
N^\perp(\gamma)/N^\beta(\gamma)$, coming from the restriction of the
projection $N^\perp(\gamma)\rightarrow
N^\perp(\gamma)/N^\beta(\gamma)$ (we denote by $N^\alpha(\gamma)$,
resp. $N^\beta(\gamma)$, the subbundle of the normal bundle of
$\gamma$, such that its fiber at $x\in\gamma$ is
$F^\alpha_x/T_x\gamma$, resp. $F^\beta_x/T_x\gamma$).

We have thus $N^\alpha(\gamma)\simeq \o1$, and all the 1-jets of
$N^\alpha(\gamma)$ yield global sections, thus local Jacobi fields.

But the existence of a Jacobi field in $\j^\alpha_\gamma\smallsetminus
\j^\gamma_\gamma$ implies, by the Jacobi equation $\ddot
J=R(\dot\gamma,J)\dot\gamma$, that $W^+(F^\alpha)=0$, the self dual
Weyl tensor vanishes on $F^\alpha$, the $\alpha$-plane generated by $\dot\gamma$.

We recall now that, for a fixed point $x$, $W^+_x$ is a polynomial of
order 4, and it is zero on $F^\alpha$ for the compact null-geodesic
$\gamma$. The same is true for the compact deformations of $\gamma$,
which implies that $W^+_x=0$. This holds for all the points of
$\gamma$, and also for all points covered by the deformations of
$\gamma$. As these deformations cover at least an open set around
$\gamma$, we conclude that $W^+$ vanishes on a non-empty open set,
thus, being holomorphic, $W^+=0$ on the whole (connected) manifold
{\bf M}.

This completes the proof of Theorem \ref{compact}
\end{proof}

\section[The projective structure of $\beta$-surfaces]{The projective structure of $\beta$-surfaces in a self-dual manifold}\label{s.w6}

The null-geodesics contained in a $\beta$-surface $\beta$ define a projective structure on
the totally-geodesic surface $\beta$, which is also given by any
connection on $\beta$ induced by a 
Levi-Civita connection on {\bf M}. We claim that this projective structure
is {\it flat}, i.e. locally equivalent to $\cp2$. 

{\noindent\bf Example. } If $\mathbf{M}=\M$, then a $\beta$-surface indexed by
$(L,l)\in\mathcal{F}$ is $\beta^{(L,l)}=\{(A,a)|A\subset l, L\subset a,
A\not\subset a\}\simeq \C^2$, and the null-geodesics in $\beta^{(L,l)}$ are
identified to the affine lines in $\C^2$ (see Section \ref{ss.w85}).
\medskip

To prove the projective flatness of a 2-dimensional manifold $\beta$,
we need to prove that the Thomas tensor $T$ vanishes identically
\cite{tho}. This tensor is an analog of the Cotton-York tensor in
conformal geometry (there is also a Weyl tensor of a projective
structure, but it only appears in dimensions greater than 2).

For a connection $\nb$ in the projective class of $\beta$, the Thomas
tensor is defined as follows \cite{tho}: For $X,Y,Z\in T\beta,$
\begin{equation}
  \label{eq:thomas}\begin{array}{rl}
  T(X,Y,Z):=&\!\!\!\!-2(\nb_ZK)(Y)X+2(\nb_YK)
  (Z)X-\\ 
&\!\!\!\!-\ \,(\nb_ZK)(X)Y\!+\
  \,(\nb_YK)(X)Z,
\end{array}
\end{equation}
where the derivation involves only the curvature tensor $R$, and the
$K(Y)X:=\mbox{tr}R(Y,\cdot)X$ is the trace of the
endomorphism $R(Y,\cdot)X\in\mbox{End}(T\beta)$.

The Thomas tensor is independent of the connection $\nb$, therefore we
will consider that $\nb$ is induced by a Levi-Civita connection
on {\bf M}.

\begin{prop}\label{p9}
  The Thomas tensor of a $\beta$-surface can be expressed in terms of the
  anti-self-dual Cotton-York tensor of {\bf M}, thus it is identically
  zero.
\end{prop}
\begin{proof}
We need first to define the {\it anti-self-dual}
  Cotton-York tensor as an irreducible component of the Cotton-York
  tensor of {\bf M}.
 
{\bf Convention} We note $C$ the Cotton-York tensor of $(\mathbf{M},c)$; we
will not use this letter for the isotropic cone in this Section, nor
in the following one.

The Cotton-York tensor is not conformally invariant; its definition
depends on a (local) metric $g$ in the conformal structure, which is
supposed to be fixed \cite{gconf}:%
\begin{equation}
  \label{eq:CY4}
  C(X,Y)(Z):=(\nb_Xh)(Y,Z)-(\nb_Yh)(X,Z),\;\forall X,Y,Z\in T\mathbf{M},
\end{equation}
where $h$ is the normalized Ricci tensor of {\bf M},
\begin{equation}
  \label{eq:rich}
  h=\frac{1}{2n(n-1)}Scal\cdot g+\frac{1}{n-2}Ric_0,
\end{equation}
$Ric_0,Scal$ being the trace-free Ricci tensor, resp. the scalar
curvature of the metric $g$, and $n:=\mbox{dim}\mathbf{M}$. In our case, $n=4$,
but the formula applies in all dimensions greater than 2 \cite{gconf}.

\obs The Cotton-York tensor $C$ of {\bf M} is a 2-form with values in $T^*\mathbf{M}$,
thus it has two components $C^+\in T^*\mathbf{M}\otimes\Lambda^+\mathbf{M}$, and $C^-\in
T^*\mathbf{M}\otimes\Lambda^-\mathbf{M}$. 
$C$  satisfies a {\it first} Bianchi identity,
as $h$ is a symmetric tensor, and also a {\it contracted}
(second) Bianchi identity, coming from the second Bianchi identity in
Riemannian geometry, \cite{gconf} :
\begin{eqnarray}\label{ar1}
\sum C(X,Y)(Z)&=&0 \;\mbox{ circular sum};\\
\label{ar2}
\sum C(X,e_i)(e_i)&=&0\;\mbox{ trace over an orthonormal basis}.
\end{eqnarray}
That means that $C\in\Lambda^2\mathbf{M}\otimes\Lambda^1\mathbf{M}$, and is orthogonal
on $\Lambda^3\mathbf{M}\subset \Lambda^2\mathbf{M}\otimes\Lambda^1\mathbf{M}$ and on
$\Lambda^1\mathbf{M}$, which is identified with the image  in
$\Lambda^2\mathbf{M}\otimes\Lambda^1\mathbf{M}$ by the metric adjoint
of the contraction (\ref{ar2}). 

Now, the Hodge operator $*:\Lambda^2\mathbf{M}\rightarrow\Lambda^2\mathbf{M}$ induces a
symmetric endomorphism of $\Lambda^2\mathbf{M}\otimes\Lambda^1\mathbf{M}$, which
maps the two above spaces isomorphically into each other. This implies
that $C^+$ and $C^-$ satisfy (\ref{ar1}) and (\ref{ar2}) (note that
these two relations are equivalent in their case).

The Cotton-York tensor is related to the Weyl tensor of {\bf M} by the
formula \cite{gconf}:
\begin{equation}
  \label{eq:dw=cy}
  \delta W=C,
\end{equation}
where $\delta:\Lambda^2\mathbf{M}\otimes\Lambda^2\mathbf{M}\rightarrow
\Lambda^2\mathbf{M}\otimes\Lambda^1\mathbf{M}$ is induced by the codifferential on the
second factor, and by the Levi-Civita connection $\nb$. Then, $C^+$
has to be the component of $\delta W$ in
$\Lambda^1\mathbf{M}\otimes\Lambda^+\mathbf{M}$, and we know that the restriction of
$W^-$ to $\Lambda^2\mathbf{M}\otimes\Lambda^+\mathbf{M}$ is identically zero. This means
that
\begin{eqnarray}
  \label{ar+}
  \delta W^+&=&C^+,\;\mbox{and also}\\
  \delta W^-&=&C^-.
\end{eqnarray}

Hence, as {\bf M} is self-dual, $C^-$ vanishes identically. 

We can prove now that the Thomas tensor of a $\beta$-surface $\beta$ is
identically zero: first we prove
\begin{equation}
  \label{eq:thomh}
  K(Y)X=\mbox{tr}|_{T\beta}R(Y,\cdot)X=h(X,Y),\; \forall X,Y\in T\beta.
\end{equation}

We recall that the {\it suspension} $h\wedge I$, viewed as an
endomorphism of $\Lambda^2\mathbf{M}$, is defined by \cite{gconf}:
\begin{equation}
  \label{eq:hI}
  (h\wedge I)(X,Y):=h(X)\wedge Y-h(Y)\wedge X,\ X,Y\in T\mathbf{M},
\end{equation}
where $h$ is identified with a symmetric endomorphism of $T\mathbf{M}$.

We have then the following decomposition of the Riemannian curvature
\cite{gconf}:
$$R=h\wedge I+W^++W^-.$$
Of course, if {\bf M} is self-dual, $W^-=0$ and $W^+(X,Y)=0$ if $X,Y\in
T\beta$ (in fact, the elements in $\Lambda^2F^\beta$, for any $\beta$-plane
$F^\beta\subset T_x\mathbf{M}$, correspond to the isotropic vectors in
$\Lambda^-\mathbf{M}$), because $W^+|_{\Lambda^-M}=0$. Then, if we choose the
basis $\{X,Y\}$ in $T\beta$, we get
\begin{eqnarray*}
K(Y)X & = & \mbox{tr}|_{T\beta}(h\wedge
I)(Y,\cdot)X =\\
        & = & \mbox{the component along $X$ of } (h\wedge I)(Y,X)X=\\
        & = & h(Y,X),
\end{eqnarray*}
which proves (\ref{eq:thomh}). The Thomas tensor of the projective
structure of $\beta$ has the following expression (see (\ref{eq:thomas})):
$$
T(X,Y,Z)=\! -3(\nb_Zh)(Y,X)+3(\nb_Yh)(Z,X)\! =\! 3C(Y,Z)(X), \, \forall
X,\! Y,\! Z\!\in\!  T\!\beta, 
$$
and, as $C^+(\cdot,\cdot)(X)$ vanishes on the anti-self-dual 2-form
$Y\wedge Z$, we conclude
\begin{equation}
  \label{eq:tcy}
T(X,Y,Z)=C^-(Y,Z)(X)=0.
\end{equation}
\end{proof}

As the flatness of the projective structure on a 2-dimensional
manifold is equivalent to the vanishing of its Thomas tensor
\cite{tho}, we get
\begin{corr}\label{prflat}
  The projective structure of the $\beta$-surfaces of a self-dual complex manifold
  {\bf M} is flat.
\end{corr}

From the classification of projectively flat compact complex surfaces
(\cite{koboch}, see also \cite{kling}), we get then a classification
of compact $\beta$-surfaces in {\bf M}: 

\begin{thrm}
  \label{clasb} A compact $\beta$-surface of a self-dual complex 4-manifold
  belongs (up to finite covering) to the following classes:
\begin{enumerate}
\item $\cp2$;
\item a compact quotient of the complex-hyperbolic plane
  $\mathbf{H}^2_\C/\Gamma$;
\item a compact complex surface admitting a (flat) affine structure: 
\end{enumerate}
(i) a Kodaira surface;\\
(ii) a properly elliptic surface with $b_1$ odd;\\
(iii) an {\it affine} Hopf surface;\\
(iv) an Inoue surface;\\
(v) a complex torus.
\end{thrm}
See \cite{kling}, \cite{koboch}, \cite{kod1} for details.




\section{Umbilic hypersurfaces in self-dual 4-manifolds}\label{s.w7}

It has been shown by LeBrun \cite{leb1} that, for a given geodesically
connected 3-manifold
with conformal structure $(\mathbf{ Q},c')$, there is a (germ-unique) self-dual
4-manifold $(\mathbf{M},c)$, such that $(\mathbf{ Q},c')$ is an umbilic 
hypersurface of $(\mathbf{M},c)$. We can therefore note by $c$ the conformal
structure of $\mathbf{ Q}$, as it coincides with the restriction of
the conformal structure of
{\bf M}.

The technical tool used in the proof of this result \cite{leb1} is the
{\it twistor space} $Z$ of $(\mathbf{ Q},c)$, which is the space of complex
null-geodesics of this manifold. It has been shown by LeBrun \cite{leb1} that
$Z$ is a 3-manifold with a {\it contact structure}, and containing
projective lines with normal bundle $\o1\oplus\o1$. 
Conversely, for any
such manifold $Z$, the space of these lines {\it tangent} to the
distribution of planes induced by the contact structure, is a
conformal 3-manifold. On the other hand, $Z$ can be
identified with the twistor space of a self-dual 4-manifold {\bf M}, in
which $\mathbf{ Q}$ is an umbilic hypersurface (the {\it conformal infinity} of
a Einstein metric on {\bf M}, \cite{leb1}). $Z$ 
has an additional structure, namely a {\it contact structure},
represented by a (non-integrable) distribution of 2-planes
$F_\beta\subset T_\beta Z$, which corresponds to the
space of Jacobi
fields $J\perp \gamma^\beta$, 
where $\gamma^\beta=\mathbf{ Q}\cap\beta$ and $Z$
is considered as the twistor space of $\mathbf{ Q}$, or, equivalently, to the
null-geodesic $\gamma^\beta\subset \mathbf{M}$, if $Z$ is 
considered as the twistor space of {\bf M}. We also remark that the above
contact structure yields a section in the bundle $B\rightarrow Z$,
which is never tangent to the $\alpha$-cones (as $\mathbf{ Q}$ is transverse to
all $\alpha$-planes ).

\obs A holomorphic contact structure on a twistor space $Z$ does not
necessarily determine a conformal 3-manifold: if $Z=\mathcal{F}$ with
its contact structure (see Section \ref{ss.w84}), there is no rational curve
in $Z$, with normal bundle $\o1\oplus\o1$, tangent to the corresponding
distribution of planes. On the other hand, this contact structure
yields an Einstein metric on $\mathbf{M}=\M$, which admits the smooth 3-manifold
$\mathcal{F}\subset \overline{\mathbf{M}}=\P(E)\times\P(E)^*$ as
``infinity''. However, this ``infinity'' is not conformal (see Section
\ref{ss.w86}).
\medskip

In the 3-dimensional case, the conformally invariant tensor
``measuring'' the non-flatness of $(\mathbf{ Q},c)$ is the Cotton-York tensor,
defined in general by 
\begin{equation}
  \label{eq:cy}
  C(X,Y):=\nabla_Xh(Y)-\nabla_Yh(X), \ X,Y\in T\mathbf{ Q},
\end{equation}
where $h$ is the {\it normalized} Ricci tensor \cite{gconf},
(\ref{eq:rich}). For the 3-manifold $\mathbf{ Q}\subset \mathbf{M}$, we have
$h(X)=\frac{1}{12}Scal+Ric_0$. 

It is natural to ask how is the Cotton-York tensor of $\mathbf{ Q}$ related to
the Weyl tensor of {\bf M}. We first recall a few facts of conformal
geometry in dimensions 3 and 4.
\smallskip

For the 3-dimensional manifold $\mathbf{ Q}$, the Riemannian curvature has the 
following expression: 
\begin{equation}\label{h3}
R^{\mathbf{ Q}}(X,Y)=(h\wedge \mbox{\bf\em I})(X,Y):= h(X)\wedge Y-h(Y)\wedge X,\;\forall X,Y\in
T\mathbf{ Q},
\end{equation}
as there is no Weyl tensor (in general, $h\wedge\mbox{\bf\em I}$
is the Ricci component of the curvature, \cite{gconf}). If we introduce
the Hodge operator $*^{\!\scriptscriptstyle{\mathbf{ Q}}}:\Lambda^2\mathbf{ Q}\rightarrow
\Lambda^1\mathbf{ Q}$, then $R^\mathbf{ Q}$ is equivalent to the symmetric 2-tensor
$*^{\!\scriptscriptstyle{\mathbf{ Q}}}\circ R^\mathbf{ Q}\circ
*^{\!\scriptscriptstyle{\mathbf{ Q}}}$. A straightforward application of the
above formula yield 
\begin{equation}\label{R*}
*R^\mathbf{ Q}:=*^{\!\scriptscriptstyle{\mathbf{ Q}}}\circ R^\mathbf{ Q} \circ
*^{\!\scriptscriptstyle{\mathbf{ Q}}}=-h+(\mathrm{tr}h)I. 
\end{equation}

For the 4-dimensional manifold {\bf M}, the components of the Riemannian 
curvature can also be expressed as eigenspaces of $*$-type operators.
Namely, regarding $R:=R^\mathbf{M}$ as a symmetric endomorphism of $\Lambda^2\mathbf{M}
=\Lambda^+\mathbf{M}\oplus\Lambda^-\mathbf{M}$, $W^+$ is the trace-free component of $R$
in $\mathrm{End}(\Lambda^+\mathbf{M})$, and $W^-$ is the trace-free component
of $R$ in $\mathrm{End}(\Lambda^-\mathbf{M})$ \cite{st}.

Let $\mathbf{ Q}\subset \mathbf{M}$ be a hypersurface, such that the restriction of the 
conformal structure $c$ of {\bf M} to $\mathbf{ Q}$ is non-degenerate (equivalently,
$T\mathbf{ Q}$ is nowhere tangent to an isotropic cone). We call $c$ the induced 
structure on $\mathbf{ Q}$. We suppose that $\mathbf{ Q}$ is umbilic.

\obs There is no possible choice of an ``orientation'' in this
case. Indeed, the group of conformal transformations of $\C^n$,
$CO(n,\C)=O(n,\C)\times \C^*/\{\pm\mathbf{1}\}$ is non-connected if
$n$ is even, and a choice of an orientation is a restriction of the
frame bundle of a $n$-dimensional conformal manifold to the connected
component of $\mathbf{1}\in CO(n,\C)$. But if $n$ is odd, $CO(n,\C)$
is connected, so all $CO(n,\C)$-frames can be connected to each other
by continuous paths. Therefore, although a complex-Riemannian
3-manifold admits, locally, 2 possible orientations, they are
conformally equivalent, fact that makes impossible a canonical way
to associate an orientation to a metric in the conformal class.

There is another way to see this difference between the even- and
odd-dimen\-sional conformal manifolds: Let $(X,c)$ be a $n$-dimensional
conformal manifold. Then a (local) metric $g$ in the conformal class $c$
is a global section in $L^2\subset S^2T^*X$, where $L$ is the bundle
of weighted scalars. We can canonically associate to $c$
a global section of $\kappa^2\otimes L^{2n}$, which is the induced metric on the
canonical bundle. For a given metric, there are only 2 (local)
sections of $\kappa$ of ``norm'' $1$. We can pick one of them if we
have a given section in $\kappa\otimes L^n$ (an ``orientation''), and
if we have a canonical way to associate to $g$,
which is a section in $L^2$, a section of $L^n$ . This can be done if
$n$ is even, namely $g^{\frac{n}{2}}$. We also remark that, if $n$ is
even, all we need to define a conformal structure on $X$ is just
the bundle $L^2$ (which is a $n/2$-th root of $\kappa$), but
if $n=2k+1$ is odd, $L^2$ automatically gives us $L\simeq\kappa\otimes
L^{\otimes -2k}$.
\medskip

We can canonically identify $\Lambda^+\mathbf{M}$ and $\Lambda^-\mathbf{M}$, restricted
to $\mathbf{ Q}\subset \mathbf{M}$, to $\Lambda^2\mathbf{ Q}$, by:
\begin{equation}
  \label{eq:star}
  \begin{array}{lcr}
\Lambda^2\mathbf{ Q}\ni\alpha& \mapsto &\alpha+*^{\!\scriptscriptstyle{\mathbf{M}}}
\alpha\in\Lambda^+\mathbf{M}\\
\Lambda^2\mathbf{ Q}\ni\alpha& \mapsto &\alpha-*^{\!\scriptscriptstyle{\mathbf{M}}}
\alpha\in\Lambda^-\mathbf{M}.
  \end{array}
\end{equation}

\begin{thrm}\label{umbilic}
  Let $\mathbf{ Q}$ be an umbilic hypersurface of a self-dual manifold
  {\bf M}. Then:

(i) The Weyl tensor of {\bf M}, restricted to $\mathbf{ Q}$, is identically zero;

(ii) The Cotton-York tensor of $\mathbf{ Q}$ is related to the self-dual Weyl
tensor of {\bf M} by the formula:
$$=g(\nb_\nu W^+(A),B)_x=-C(A)(*^{\!\scriptscriptstyle{\mathbf{ Q}}}B)_x,$$
where $A,B\in\Lambda^2_x\mathbf{ Q}$, $\nu\perp T_x\mathbf{ Q}$ is unitary for the metric
$g$, and the Hodge operator $ *^{\!\scriptscriptstyle{\mathbf{ Q}}}$ is induced
by $g$ and the orientation on $\mathbf{ Q}$ admitting $\nu$ as an exterior
normal vector.
\end{thrm}
\begin{corr}\label{cumb}
  If $(\mathbf{M},c)$ is self-dual and $\mathbf{ Q}\subset \mathbf{M}$
  is an umbilic hypersurface, 
  then the Cotton-York tensor of $\mathbf{ Q}$, $C^\mathbf{ Q}$, is identified to the
  restriction to $\mathbf{ Q}$ of the
  self-dual Cotton-York tensor $C^+$ of {\bf M}:
$$C^+(X,Y)(Z)=C^\mathbf{ Q}(X,Y)(Z),\;\forall X,Y,Z\in T\mathbf{ Q}.$$
\end{corr}
{\noindent\sc Proof of the Theorem. } The claimed identity is conformally invariant: If
$X,Y,Z,\nu$ is a $g$-ortho\-normal oriented basis of {\bf M}, then $X,Y,Z$ is
a $g$-orthonormal basis on $\mathbf{ Q}$ giving the orientation as above. Then
$*^{\!\scriptscriptstyle{\mathbf{ Q}}} (Z\wedge X)=Y$, and, if we take
$A:=X\wedge Y,\, B:=Z\wedge X$, the claimed identity becomes 
\begin{equation}
  \label{eq:CYW}
  \langle \nb_\nu W^+(X,Y)Z,X\rangle =-C(X,Y)(Y),
\end{equation}
where angle brackets denote the scalar product induced by $g$.
 
The tensors $W^+, C$, in the above form, are independent of the chosen
metric $g$ \cite{gconf}, which depends on the normal vector $\nu$,
supposed to be $g$-unitary. If $\nu':=\lambda\nu$, for
$\lambda\in\C^*$, then the corresponding metric $g'=\lambda^{-2}g$,
and also ${*^{\!\scriptscriptstyle{\mathbf{ Q}}}}'=\lambda^{-1}
*^{\!\scriptscriptstyle{\mathbf{ Q}}}$, thus the identity (\ref{eq:CYW}) for
$\nu',g'$ is equivalent to the one for $\nu,g$.

\obs As $W^+$ is the trace free component of the Riemannian curvature
contained in $\mbox{End}(\Lambda^+\mathbf{M})$, and is symmetric, it is enough
to evaluate it on pairs $A,B\in\Lambda^2\mathbf{ Q}\simeq\Lambda^+\mathbf{M}$ which are
unitary and orthogonal for the metric $g$, therefore the 
check of the equation (\ref{eq:CYW}) will prove the theorem.

As $W^\pm$ are $*^{\!\scriptscriptstyle{\mathbf{M}}}$-eigenvectors in
$\mbox{End}_0(\Lambda^2\mathbf{M})$ (the space of trace-free endomorphisms of
$\Lambda^2\mathbf{M}$), they are determined by the following formulas, where
$X,Y,Z$ is any oriented orthonormal basis of $T\mathbf{ Q}$:
\begin{eqnarray}
  \label{arw+}
 \  \langle W^+(X,Y)Z,X\rangle &\!\! =&\!\!\frac{1}{4}(\langle R(X,Y)Z,X\rangle +\langle R(Z,\nu)Y,\nu\rangle +\\
\nonumber  &&+\langle R(X,Y)Y,\nu\rangle + \langle R(Z,\nu)Z,X\rangle )\\
\label{arw-}
 \  \langle W^-(X,Y)Z,X\rangle &\!\! =&\!\!\frac{1}{4}(\langle R(X,Y)Z,X\rangle +\langle R(Z,\nu)Y,\nu\rangle -\\
\nonumber  &&-\langle R(X,Y)Y,\nu\rangle - \langle R(Z,\nu)Z,X\rangle ), 
\end{eqnarray}
where $X,Y,Z,\nu$ is supposed to be a local extension, around a region
of $\mathbf{ Q}$, of the $g$-orthonormal frame used in (\ref{eq:CYW}). As {\bf M} is
supposed to be self-dual, $W^-$ is identically zero, thus, in the
points $x\in \mathbf{ Q}$, we have  
\begin{equation}
  \label{eq:W1}
  \langle W^+(X,Y)Z,X\rangle _x=\frac{1}{2}(\langle R(X,Y)Y,\nu\rangle + \langle R(Z,X)Z,\nu\rangle )_x.
\end{equation}

It is a standard fact that, if $\mathbf{ Q}$ is umbilic, there is a local metric
$g$ in the conformal class $c$ of {\bf M}, such that, for $g$, $\mathbf{ Q}$ is
totally geodesic. We fix such a metric. Then we have 
\begin{equation}
  \label{eq:tgeod}
  R(X',Y')Z'=R^\mathbf{ Q}(X',Y')Z',\;\forall X',Y',Z'\in T\mathbf{ Q},
\end{equation}
which, together with \ref{eq:W1}, implies that $W^+|_\mathbf{ Q}\equiv 0$.
\smallskip

On the other hand, (\ref{eq:tgeod}), together with (\ref{eq:W1}) and
(\ref{arw-}), yield
\begin{equation}
  \label{eq:rq}
  \langle R(X,Y)Z,X\rangle _x+\langle R(Z,\nu)Y,\nu\rangle _x=0,
  \forall x\in{\bf Q}.
\end{equation}

Let us compute now the normal derivative of $W^+$ in a point $x\in \mathbf{ Q}$;
we suppose that $X,Y,Z,\nu$ are locally extended by an orthonormal
frame, and that they are parallel at $x$ (we omit, for simplicity of
notation, the point $x$ in the following lines:
$$
\langle \nb_\nu W^+(X,Y)Z,X\rangle  =\frac{1}{2}(\langle \nb_\nu R(X,Y)Y,\nu\rangle +\langle \nb_\nu
R(Z,X)Z,\nu\rangle ),$$
from (\ref{eq:W1}). This is then equal to:
$$\begin{array}{rcl}
\langle \nb_\nu W^+(X,Y)Z,X\rangle &\!\!\! =&\!\!\!-\frac{1}{2}(\langle \nb_X
R(Y,\nu)Y,\nu\rangle +\langle \nb_Y
R(\nu,X)Y,\nu\rangle +\\
 &&+\langle \nb_Z R(X,\nu)Z,\nu\rangle  +\langle \nb_X R(\nu,Z)Z,\nu
 \rangle ),\end{array}$$  
from the second Bianchi identity. Then we have
$$\begin{array}{rcl}
\langle \nb_\nu W^+(X,Y)Z,X\rangle &\!\!\! =&
\!\!\!\frac{1}{2}(\langle \nb_X R(Z,X)Z,X\rangle +\langle \nb_Y R(Z,Y)Z,X\rangle +\\
&&+\langle \nb_Z R(Y,Z)X,Y\rangle +\langle \nb_X R(Y,X)X,Y\rangle )\end{array}$$
from analogs of (\ref{eq:rq}). Then
$$\begin{array}{rcl}
\langle \nb_\nu W^+(X,Y)Z,X\rangle &\!\!\!\!\! =&\!\!\!\!\frac{1}{2}(\langle \nb_X
R^\mathbf{ Q}(Z,X)Z,X\rangle  +\langle \nb_Y R^\mathbf{ Q}(Z,Y)Z,X\rangle  +\\
&&\!\!\! +\langle \nb_Z R^\mathbf{ Q}(Y,Z)X,Y\rangle  +\langle \nb_X
R^\mathbf{ Q}(Y,X)X,Y\rangle ) 
\end{array}$$
from (\ref{eq:tgeod})
$$\begin{array}{rcl}
\langle \nb_\nu W^+(X,Y)Z,X\rangle &\!\!\! =&\!\!\!\frac{1}{2}(\nb_X
h(Z,Z)+\nb_X h(X,X)+\nb_Y h(Y,X)-\\ 
&&-\nb_Z h(Z,X)-\nb_X h(X,X)-\nb_X h(Y,Y)),
\end{array}$$
from (\ref{h3}). Finally, from (\ref{ar2}), we get
$$
\langle \nb_\nu W^+(X,Y)Z,X\rangle =\frac{1}{2}(C(X,Z)(Z)-C(X,Y)(Y))=-C(X,Y)(Y)$$
This proves equation (\ref{eq:CYW}).
\endproof

The Corollary \ref{cumb} now easily follows from the above theorem and
(\ref{ar+}). 
\bigskip 

The main results in Section \ref{s.w5} hold also in the case of a
conformal (complex) 3-manifold:

\begin{thrm}\label{compact3}
  Let $Z$ be a twistor space of a conformal 3-manifold $\bf Q$; let
  $F_{\bar\gamma}\subset T_{\bar\gamma}Z$ be its contact
  structure. Suppose there is a point $\bar\gamma\in Z$ such that,
  to any direction in $F_{\bar\gamma}$, there is a tangent rational
  curve in $Z$ with normal bundle $\o1\oplus\o1$. Then $Z$ is
  projectively flat, and $\mathbf{ Q}$ is conformally flat.
\end{thrm}
This follows directly from Theorem \ref{twist}, as $\mathbf{ Q}$ is umbilic in
{\bf M}, the space of all projective lines in $Z$ with the above normal
bundle \cite{leb1}, and from Theorem \ref{umbilic}.

\begin{thrm}
  \label{comp3} Let $\mathbf{ Q}$ be a conformal 3-manifold containing an
  immersed rational curve as null-geodesic. Then $\mathbf{ Q}$ is conformally flat.
\end{thrm}

\begin{proof} We cannot use directly Theorem \ref{umbilic}, as the ``ambient''
self-dual manifold {\bf M} can only be defined for a {\it civilized}
(e.g. geodesically connected) 3-manifold. Hence, we follow the steps
in Proposition \ref{tub} and prove 

\begin{prop}\label{tub3}
Null-geodesics close to a compact, simply-connected one are also
compact and simply-connected, and they are embedded.
\end{prop}

Still using the same arguments as in Section \ref{s.w5}, we get an embedded null-geodesic
$\gamma\subset \mathbf{ Q}$ diffeomorphic to $\cp1$, and we have:

\begin{prop}
  The deformations of $\gamma$ as a compact curve coincide with the
  null-geodesics close to $\gamma$.
\end{prop}

We cover $\gamma$ with geodesically convex open sets $U_i, \,
i=\overline{1,n}$, such that:
\begin{equation}
  \label{eq:gcx3}
  \forall i\neq j \mbox{ s.t. }U_i\cap U_j\cap\gamma\neq\emptyset,
  \exists U_{ij}\supset(U_i\cup U_j),
\end{equation}
where $U_{ij}$ is still geodesically convex (with respect to a
particular Levi-Civita connection). This is possible by choosing $U_i,
\, i=\overline{1,n}$, small enough. Then we choose a relatively
compact tubular neighborhood $N(r_0)$ of $\gamma$, such that its
closure is covered by the $U_i$'s. 

We consider then the twistor spaces $Z_i$, the spaces of null-geodesics of
$U_i$. The compact, simply-connected, null-geodesics close to $\gamma$ identify
(diffeomorphically) the neighborhoods of $\bar\gamma^i\in Z_i$ with
the space $Z$ of the deformations of $\gamma$ as a compact curve. We
can see then $Z$ as an open set common to all the $Z_i$'s:
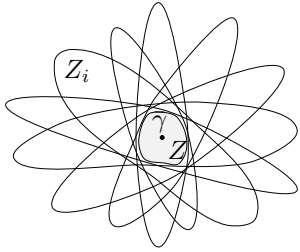
\begin{center}\begin{picture}(0,0)%
\epsfig{file=weyl/3tw.pstex}%
\end{picture}%
\setlength{\unitlength}{0.00022700in}%
\begingroup\makeatletter\ifx\SetFigFont\undefined
\def\x#1#2#3#4#5#6#7\relax{\def\x{#1#2#3#4#5#6}}%
\expandafter\x\fmtname xxxxxx\relax \def\y{splain}%
\ifx\x\y   
\gdef\SetFigFont#1#2#3{%
  \ifnum #1<17\tiny\else \ifnum #1<20\small\else
  \ifnum #1<24\normalsize\else \ifnum #1<29\large\else
  \ifnum #1<34\Large\else \ifnum #1<41\LARGE\else
     \huge\fi\fi\fi\fi\fi\fi
  \csname #3\endcsname}%
\else
\gdef\SetFigFont#1#2#3{\begingroup
  \count@#1\relax \ifnum 25<\count@\count@25\fi
  \def\x{\endgroup\@setsize\SetFigFont{#2pt}}%
  \expandafter\x
    \csname \romannumeral\the\count@ pt\expandafter\endcsname
    \csname @\romannumeral\the\count@ pt\endcsname
  \csname #3\endcsname}%
\fi
\fi\endgroup
\begin{picture}(6819,5657)(2102,-6993)
\put(3496,-3076){\makebox(0,0)[lb]{\smash{\SetFigFont{10}{12.0}{rm}
$Z_i$
}}}
\put(5515,-4189){\makebox(0,0)[lb]{\smash{\SetFigFont{10}{12.0}{rm}
$\bar\gamma$
}}}
\put(5881,-4921){\makebox(0,0)[lb]{\smash{\SetFigFont{10}{12.0}{rm}
$Z$
}}}
\end{picture}
\end{center}
Following LeBrun, we define the self-dual manifolds $\mathbf{M}_i$ as the space
of projective lines in $Z_i$, with normal bundle $\o1\oplus\o1$. Then
$U_i$ is an umbilic hypersurface in $\mathbf{M}_i$. 

The local twistor spaces $Z_i$ admit contact structures, which
coincide on $Z$, and contain projective lines $Z^i_x$ corresponding to
points $x\in\gamma\cap U_i$. If we denote by $Z_{ij}$ the twistor
space of $U_{ij}$, then $Z_i$ and $Z_j$ are identified to open sets in
$Z_{ij}$, in particular the lines $Z^i_x$ and $Z^j_x$ are identified, thus their
intersections with the common set $Z$ coincide. We denote by $Z_x$ this
(non-compact) curve in $Z$, and by $F$ the canonical contact structure
of $Z$ (restricted from the ones of $Z_i$).

\obs We already have obtained that the $\alpha$-cone corresponding to
$F_{\bar\gamma}$ is a part of a smooth surface: 
the union of the lines
$Z_x$, $x\in\gamma$, thus, from Theorem \ref{3'}, the Weyl tensor $W_i^+$ of
the self-dual manifold $\mathbf{M}_i$ vanishes on the $\alpha$-planes generated by
$T\gamma$. But this is nothing new: we know, from Theorem
\ref{umbilic}, that $W_i^+$ vanishes on $U_i$.
\smallskip

We intend to apply Theorem \ref{3'} to prove that $W_i^+$ vanishes on
points close to $U_i$, but in $\mathbf{M}_i\smallsetminus U_i$. We do that by
showing that the integral $\alpha$-cones corresponding to planes $F^y\subset
T_{\bar\gamma}Z$ are parts of smooth surfaces, then we conclude using
Theorem \ref{3'}.
\smallskip

First we choose hermitian metrics $h_i$ on $Z_i$, such that they
coincide (with $h$) on $Z$. We have a diffeomorphism between $\gamma$
and $\P(F_{\bar\gamma})$, so we choose relatively compact open sets in
$\P(T_{\bar\gamma}Z)$, covering $\P(F_{\bar\gamma})$, with the
following properties: As the metrics $h_i$ induce metrics on $\mathbf{M}_i$, we
first choose a small enough distance $r_1>0$ such that 
\begin{enumerate}
\item $\forall i$, there is a sub-covering $V_i\Subset U_i$ of $\gamma$
  such that the
  ``tubular neighborhoods'' $W_i:=\{y\in \mathbf{M}_i|\mathrm{d}(y,\bar V_i)\leq
  r_1\,\, \pi_i(y)\in \bar V_i\cap\gamma\}$ are compact ($\pi_i$ is the
  ``orthogonal projection'' --- for the hermitian metric --- from
  $\mathbf{M}_i$ to $\gamma\cap U_i$; it is well defined because of the
  condition below); 
\item $r_1$ is less than the bijectivity radius of the (hermitian)
  exponentials for the points of $\bar V_i$ in $\mathbf{M}_i$, and for the
  points of $\overline{V_i\cup V_j}$ in $\mathbf{M}_{ij}$ (if $U_i\cap
  U_j\cap\gamma\neq\emptyset.$). 
\end{enumerate}
We have then
\begin{lem} For any $y_i\in W_i\subset \mathbf{M}_i$, $y_j\in W_j\subset \mathbf{M}_j$
  such that the curves $Z_{y_i}:=Z^i_{y_i}\cap Z$,
  $Z_{y_j}:=Z^j_{y_j}\cap Z$ are tangent to the same direction in
  $\bar\gamma\in Z$, the respective curves $Z_{y_i}, Z_{y_j}$
  coincide.
\end{lem}
\begin{proof}
We first note that the projection $\pi_i$ from $\mathbf{M}_i$ is equivalent to
the $h$- orthogonal projection of the direction of
$T_{\bar\gamma}Z_{y_i}$ to a direction in $F_{\bar\gamma}$, so
$\pi_i(y_i)=\pi_j(y_j)=:y\in\gamma$; thus $y$ belongs to both $U_i$
and $U_j$, and we use again the twistor space $Z_{ij}$ to conclude
that $Z_{y_i}$ and $Z_{y_j}$ are ``restrictions'' to $Z$ of the same
projective line (as they both have the same tangent space at
$\bar\gamma$) $Z^{ij}_{y_{ij}}$, for a point $y_{ij}\in \mathbf{M}_{ij}$.
\end{proof} 
 
Now we have a tubular neighborhood $S\subset \P(T_{\bar\gamma}Z)$ of
$\P(F_{\bar\gamma})$, of radius $r_1/2$, such that, for any 2-plane
$F'\subset S$, the conditions in Theorem \ref{3'} are satisfied
(considering any of the local twistor spaces $Z_i$). 

We conclude that the Weyl tensor $W_i^+$ of $\mathbf{M}_i$ vanishes along all
null-geodesics of $\mathbf{M}_i$, close (in $W_i$) to $\gamma$ and included in the $\beta$-surface
$\beta^i$, determined by $\gamma$. This means that $W^+$ vanishes
everywhere on $\beta^i$. By deforming $\gamma$, we obtain that $W^+_i$
vanishes on a neighborhood of $U_i$ in $\mathbf{M}_i$, hence $\mathbf{M}_i$, as well as
$U_i$, are conformally flat (by Theorem \ref{umbilic}). 

It follows from Theorem \ref{umbilic} that $\mathbf{ Q}$ is conformally flat.
\end{proof}

\section{Examples}\label{s.w8}
\subsection{The flat case }\label{ss.w81} The first example is the ``flat'' case:
$Z=\cp3=\P(\C^4)$, with its canonical projective structure, and its space of
projective lines $\mathbf{M}=\mbox{Gr}(2,\C^4)$. ($Z$ is equally the twistor
space of the Riemannian round 4-sphere, which is, therefore, a real
part of $\mbox{Gr}(2,\C^4)$.) If $\beta\in Z$, then the $\beta$-surface
associated to it is the set
$\{x\in\mbox{Gr}(2,\C^4)|\beta\subset x\subset\C^4\}$. In this flat case, we can
equally define the $\alpha$-twistor space $Z^*$, which is the dual
projective 3-space $(\cp3)^*:=\P((\C^4)^*)=\mbox{Gr}(3,\C^4)$, and an
$\alpha$-surface $\alpha\in Z^*$ is the set
$\{x\in\mbox{Gr}(2,\C^4)|x\subset\alpha\subset\C^4\}\subset \mathbf{M}$. A null-geodesic
$\gamma$ is then determined by a pair of {\it incident} $\alpha$-,
resp. $\beta$-surface $\beta\subset\alpha\subset\C^4$:
$$\gamma=\{x\in\mbox{Gr}(2,\C^4)|\beta\subset x\subset\alpha\}.$$
$\alpha$-surfaces and $\beta$-surfaces are diffeomorphic to $\cp2$, null-geodesics to $\cp1$, and the
ambitwistor space $B$ is the ``partial flag'' manifold
$$B=\{(\alpha,\beta)\in(\cp3)^*\times\cp3| \beta\subset\alpha\}.$$
The flag manifold, of dimension 7, is isomorphic to total space of the
projective cone bundle over {\bf M}, $\P(C)$. 
\smallskip

\subsection{$\cp2$ }\label{ss.w82} Another example is when $Z$ is the
twistor space of the 
real Riemannian manifold $\cp2$, with the Fubini-Study metric. Then
$Z$ is the manifold
of flags in $E=\C^3$, $\mathcal{F}:=\{(L,l)\in\P(E)\times\P(E)^*|L\subset l\}$,  ($\P(E)$,
resp. $\P(E)^*$ are viewed as the space of lines, resp. 2-planes, in
$E$) \cite{ahs}. A projective line $Z_x$ in $Z$ is a set 
$$Z_x=\{(L,l)\in\mathcal{F}|L\subset a^x,\; A^x\subset l\},
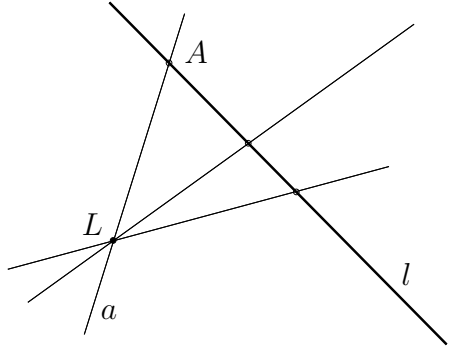
$$
\begin{center}\begin{picture}(0,0)%
\epsfig{file=weyl/cp2.pstex}%
\end{picture}%
\setlength{\unitlength}{0.00043700in}%
\begingroup\makeatletter\ifx\SetFigFont\undefined
\def\x#1#2#3#4#5#6#7\relax{\def\x{#1#2#3#4#5#6}}%
\expandafter\x\fmtname xxxxxx\relax \def\y{splain}%
\ifx\x\y   
\gdef\SetFigFont#1#2#3{%
  \ifnum #1<17\tiny\else \ifnum #1<20\small\else
  \ifnum #1<24\normalsize\else \ifnum #1<29\large\else
  \ifnum #1<34\Large\else \ifnum #1<41\LARGE\else
     \huge\fi\fi\fi\fi\fi\fi
  \csname #3\endcsname}%
\else
\gdef\SetFigFont#1#2#3{\begingroup
  \count@#1\relax \ifnum 25<\count@\count@25\fi
  \def\x{\endgroup\@setsize\SetFigFont{#2pt}}%
  \expandafter\x
    \csname \romannumeral\the\count@ pt\expandafter\endcsname
    \csname @\romannumeral\the\count@ pt\endcsname
  \csname #3\endcsname}%
\fi
\fi\endgroup
\begin{picture}(5295,4161)(2104,-4939)
\put(3001,-3614){\makebox(0,0)[lb]{\smash{\SetFigFont{12}{14.4}{rm}$L$}}}
\put(6781,-4199){\makebox(0,0)[lb]{\smash{\SetFigFont{12}{14.4}{rm}$l$}}}
\put(3226,-4604){\makebox(0,0)[lb]{\smash{\SetFigFont{12}{14.4}{rm}$a$}}}
\put(4216,-1574){\makebox(0,0)[lb]{\smash{\SetFigFont{12}{14.4}{rm}$A$}}}
\end{picture}
\end{center}
where $(A^x,a^x)$ belongs to $\P(E)\times\P(E)^*\smallsetminus \mathcal{F}$,
which is, therefore, the space {\bf M} of such lines, and a conformal
self-dual 4-manifold. It can be naturally compactified within the
space of analytic cycles of $Z$ to $\overline{\mathbf{M}}=\P(E)\times\P(E)^*$, which is
obviously a smooth manifold, but it carries no global conformal structure, as
its canonical bundle has no square root. This means that the conformal
structure on $\overline{\mathbf{M}}$ is smooth on {\bf M}, and {\it singular} on
$\mathcal{F}=\overline{\mathbf{M}}\smallsetminus \mathbf{M}$. The cycles of $Z$ corresponding to a point
$\bar x=(A,a)$ in this subset are pairs of complex projective lines in $Z$:
$$Z_{\bar x}=\{(A,l)\in Z=\mathcal{F}\}\cup\{(L,a)\in
Z=\mathcal{F}\}.$$
A $\beta$-surface in {\bf M}, corresponding to a point $\beta=(L,l)\in
Z$, is the set $$\beta=\{(A,a)\in\P(E)\times\P(E)^*|A\subset l,\;
L\subset a,\; A\ne L,\; a\ne l\},$$
and can be naturally compactified to 
$$\bar\beta=\{(A,l^\beta)\in\mathcal{F}\}\times
\{(L^\beta,a)\in\mathcal{F}\}\simeq\cp1\times\cp1.$$

\subsection{The tangent space to $\mathcal{F}$}\label{ss.w83} In order to describe
the null-geodesics of {\bf M} as 2-planes in $Z$, we study
first the tangent space of $Z=\mathcal{F}$ at $\beta=(L,l)$:

A vector in $T_{(L,l)}\mathcal{F}$ is a pair of vectors $(V,v)$, with
$V\in T_L\P(E)$ and $v\in T_l\P(E)^*$, which satisfy a linear
condition (as $\mathcal{F}\subset\P(E)\times\P(E)^*$). Actually, there
is a duality between $\P(E)^*$, the Grassmannian of 2-planes in $E$,
and $\P(E^*)$, the projective space of $E^*:=\mbox{Hom}(E,\C)$, and an
analogous one between $\P(E)$ and $\P(E^*)^*$:
\begin{eqnarray*}\P(E)^*\ni l&
  \stackrel{\simeq}{\longmapsto}&l^o\in\P(E^*) \\
\P(E)\ni L& \stackrel{\simeq}{\longmapsto}&L^o\in\P(E^*)^*.
\end{eqnarray*}
Then, the flag manifold $\mathcal{F}$ is defined, as a submanifold of
$\P(E)\times\P(E)^*$, by the equation
$$ y(Y)=0,\ \forall y\in l^o,\forall Y\in L.$$
The vector $V\in T_L \P(E)$ is a homomorphism in
$\mbox{Hom}(L,E/L)$. By duality, $v\simeq v^0\in \mbox{Hom}(l^o,E^*/l^o)$.
Then the vector $(V,v)\in T_{(L,l)}\P(E)\times\P(E)^*$ lies in
$\mathcal{F}$ iff:
\begin{equation}
  \label{eq:tflag}
  v^o(y^o)(Y)+y^o(V(Y))=0,\ \forall Y\in L,\,\forall y^o\in l^o,
\end{equation}
or, equivalently, 
\begin{equation}
  \label{eq:tflag1}
v|_L=\pi_l\circ V,
\end{equation}
where $\pi_l:E/L\rightarrow E/l$ is the projection (as $L\subset l$).

The geometry of $\mathcal{F}$, as a subset of $\P(E)\times\P(E)^*$,
can be described in the following figure:
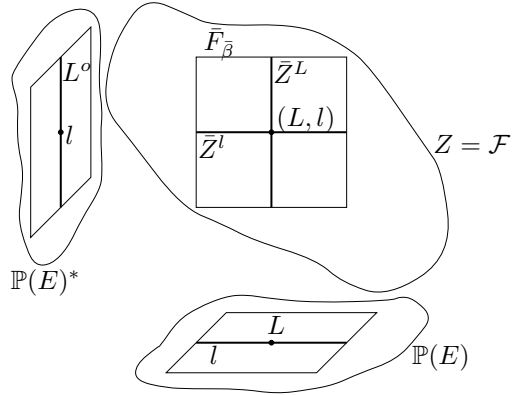
\begin{center}\input{weyl/twcp2.pstex_t}\end{center}
\subsection{The 2-planes in $\mathcal{F}$}\label{ss.w84}
Let us consider now a 2-plane $F$ in $T_{(L,l)}\mathcal{F}$, and the
cycles (corresponding to points in $\overline{\mathbf{M}}$) tangent to it. We have three
cases: \smallskip

{\noindent\bf 1. } $F=\bar F_\beta$ is the ``degenerate'' 2-plane
tangent to the 2 
special curves $\bar Z_L$, $\bar Z_l$ whose union is the special cycle
$\bar F_(L,l)$ corresponding to $(L,l)\in\overline{\mathbf{M}}\smallsetminus \mathbf{M}$. There are
no projective lines $Z_x,\, x\in \mathbf{M}$, tangent to it; only the special
cycles $\bar Z_{(L,a)},\, L\subset a$ and $\bar Z_{(A,l)},\, A\subset
l$ are tangent to $\bar F_{(L,l)}$, actually only to the two
privileged directions of $\bar Z_L$, resp. $\bar Z_l$.

\obs The special curves $\bar Z_L$, $\bar Z_l$ have trivial normal
bundle, being fibers of the projections from $\mathcal{F}$ to $\P(E)$,
resp. $\P(E)^*$, so these special curves form two complete families of
analytic cycles in $\mathcal{F}$, isomorphic to $\P(E)$,
resp. $\P(E)^*$. Two such curves are incident iff they are of
different types ($\bar Z_L$ is {\it of type} $E$, $\bar Z_l$ is {\it
  of type} $E^*$), so they can only form ``polygons'' with an even
number of edges. But there are no quadrilaterals, as one can easily
check, using the fact that $\bar Z_L$ and $\bar Z_l$ are incident iff
$L\subset l$, thus iff $l$ is a line in $\P(E)$ containing $L$. But
there are hexagons, corresponding to the 3 vertices and 3 sides of a
triangle in $\P(E)\simeq\cp2$:
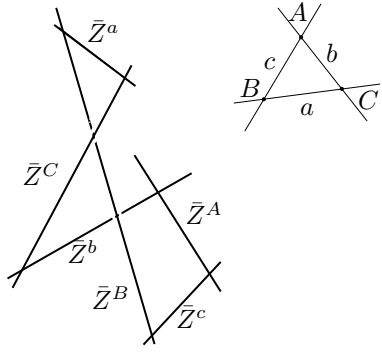
\begin{center}\input{weyl/hex1.pstex_t}\end{center}
The hexagon above is not ``flat'', i.e. there is no canonical
submanifold of $\mathcal{F}$ containing it. This, and the fact that
there are no quadrilaterals made of $\bar Z$-type curves, is just a
consequence of the fact that the distribution $\bar F$ on
$Z=\mathcal{F}$ is {\it non integrable}; in fact it is the holomorphic
{\it contact structure} induced by the Fubini-Study {\it Einstein
  metric} on $\cp2$, \cite{besse}, see also Section \ref{ss.w86}.
\smallskip

{\noindent\bf 2. } $F=F^a$, for $a\supset L$, $a\neq l$. This is a
2-plane that is tangent to only one of the special curves $\bar
Z_L$. The projective lines tangent to $F^a$ at $\beta=(L,l)$ are
$Z_{(A,a)}$, $\forall A\subset l,\, A\neq L$, hence the corresponding
null-geodesic is 
\begin{equation}\label{geodspes}
\gamma^a=\{(A,a)\in\P(E)\times\P(E)^*|A\subset l,\, A\neq l\},
\end{equation}
thus it is diffeomorphic to $\C$, and its closure is 
$$\bar
\gamma^a=\{(A,a)\in\P(E)\times\P(E)^*|A\subset l\}\simeq\cp1.$$ 
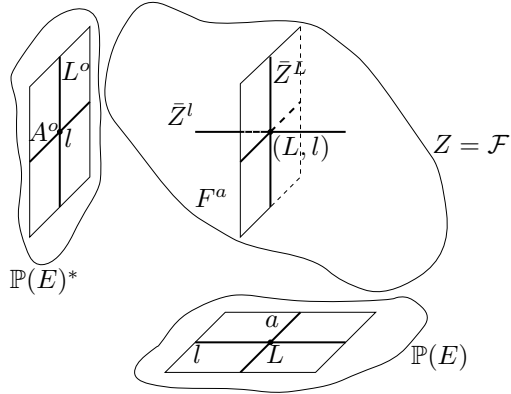
\begin{center}\input{weyl/geodspec.pstex_t}\end{center}
\obs The ``limit'' curve is $\bar Z_{(L,a)}$, so it is
non-singular at $(L,l)$. Actually, the points of $Z_{(A,a)}$ {\it
  close} to $(L,l)$ converge, when $A\rightarrow L$, to some points in
$\bar Z_L$, which is tangent to $F^a$. We can, then, apply the same
method as in Theorem \ref{twist} to conclude that the integral
$\alpha$-cone associated to $F^a$ is a smooth manifold {\it around}
$(L,l)$, thus, from Theorem \ref{t1}, the Weyl tensor $W^+$ of {\bf M}
vanishes on the $\beta$-planes generated, along $\gamma^a$, by its own
direction. We will see that the vanishing of $W^+$ on these $\alpha$-planes leads
to the existence of some $\alpha$-surfaces, see below. Of course, the deformation
argument in Theorem \ref{twist} does not hold in the present case, as
the normal bundle of $\bar Z_L$ is trivial, thus different from the
one of the rest of the rational curves $Z_{(A,a)}$ (as we will see below,
generic 2-planes through $(L,l)$ do not admit projective lines tangent
to all their directions).
\smallskip

{\noindent$\bf 2'.$ } We have a similar situation for planes $F=F^A$ ---
$A\subset l$, $A\neq L$ --- tangent to the other special curve $\bar
Z_l$.
\smallskip

{\noindent\bf 3. } $F=F^\varphi$, where $\varphi
:\P(l)\rightarrow\P(L^o)$ is a projective diffeomorphism such that
$\varphi(L)=l^o$. Indeed, the
tangent spaces $T_L\P(E)$ and $T_l\P(E)^*$ are isomorphic to
$\mbox{Hom}(L^o,E^*/L^o)$, resp. to $\mbox{Hom}(l,E/l)$, and a generic
2-plane $F$ in $T_{(L,l)}\mathcal{F}$ is the graph of a linear
isomorphism $\phi:T_L\P(E)\rightarrow T_l\P(E)^*$ satisfying a linear
condition (\ref{eq:tflag}) or (\ref{eq:tflag1}). Actually, the graph is
determined by the projective application $\varphi$ induced by $\phi$ from
$\P(T_L\P(E))\simeq \P(L^o)$ to $\P(T_l\P(E)^*)\simeq\P(l)$:
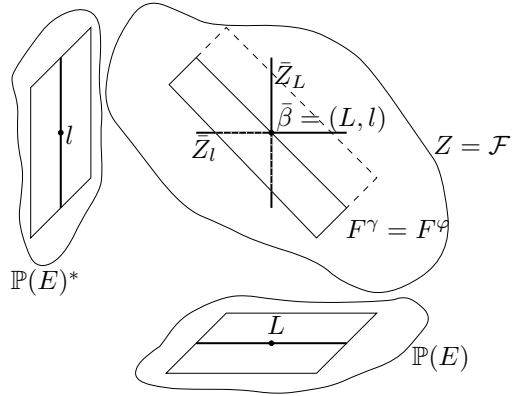
\begin{center}\input{weyl/gengeod.pstex_t}\end{center}
The condition $\varphi(L)=l^o$ is implied by (\ref{eq:tflag1}). The null-geodesic
associated to the 2-plane $F^\varphi$ is  
\begin{equation}
  \label{eq:gengeod}
  \gamma^\varphi=\{(A,a)\in\P(E)\times\P(E)^*\smallsetminus\mathcal{F}
  |A\subset l,\,a^o\subset L^o\, a^o=\varphi(A)\},
\end{equation}
and its closure in $\overline{\mathbf{M}}$ is 
\begin{equation}
  \label{eq:gengeod1}
  \bar\gamma^\varphi=\{(A,a)\in\P(E)\times\P(E)^*|A\subset l,
  \,a^o\subset L^o\},
\end{equation}
hence the ``limit'' point is $(L,l)\in\overline{\mathbf{M}}$, corresponding to the
special cycle $\bar Z_{(L,l)}$, none of whose components is tangent to
$F^\varphi$. The integral $\alpha$-cone associated
to $F^\varphi$ looks like suggested in the picture below:
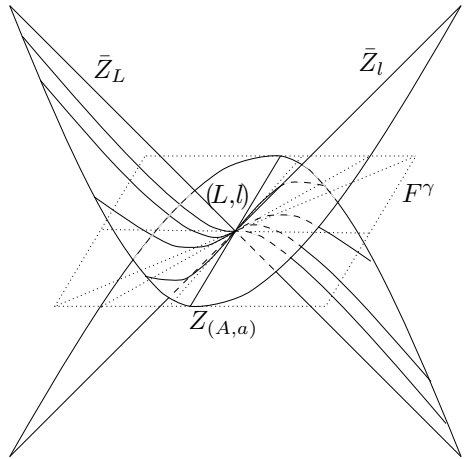
\begin{center}\begin{picture}(0,0)%
\epsfig{file=weyl/alfint.pstex}%
\end{picture}%
\setlength{\unitlength}{0.00035000in}%
\begingroup\makeatletter\ifx\SetFigFont\undefined
\def\x#1#2#3#4#5#6#7\relax{\def\x{#1#2#3#4#5#6}}%
\expandafter\x\fmtname xxxxxx\relax \def\y{splain}%
\ifx\x\y   
\gdef\SetFigFont#1#2#3{%
  \ifnum #1<17\tiny\else \ifnum #1<20\small\else
  \ifnum #1<24\normalsize\else \ifnum #1<29\large\else
  \ifnum #1<34\Large\else \ifnum #1<41\LARGE\else
     \huge\fi\fi\fi\fi\fi\fi
  \csname #3\endcsname}%
\else
\gdef\SetFigFont#1#2#3{\begingroup
  \count@#1\relax \ifnum 25<\count@\count@25\fi
  \def\x{\endgroup\@setsize\SetFigFont{#2pt}}%
  \expandafter\x
    \csname \romannumeral\the\count@ pt\expandafter\endcsname
    \csname @\romannumeral\the\count@ pt\endcsname
  \csname #3\endcsname}%
\fi
\fi\endgroup
\begin{picture}(6794,6794)(654,-6383)
\put(1936,-706){\makebox(0,0)[lb]{\smash{\SetFigFont{10}{12.0}{rm}
$\bar Z_L$
}}}
\put(5896,-601){\makebox(0,0)[lb]{\smash{\SetFigFont{10}{12.0}{rm}
$\bar Z_l$
}}}
\put(3361,-4456){\makebox(0,0)[lb]{\smash{\SetFigFont{10}{12.0}{rm}
$Z_{(A,a)}$
}}}
\put(6541,-2551){\makebox(0,0)[lb]{\smash{\SetFigFont{10}{12.0}{rm}
$F^\gamma$
}}}
\put(3601,-2566){\makebox(0,0)[lb]{\smash{\SetFigFont{10}{12.0}{rm}
$(\! L,\! l\! )$
}}}
\end{picture}
\end{center}
\subsection{The null-geodesics of the complexification of $\cp2$}\label{ss.w85} 
The application $\varphi$ has the following interpretation in terms of 
projective geometry on $\cp2=\P(E)$: a direction $\C v$ in $T_l\P(E)^*$
is identified to the point $\ker v\equiv A\in l/L\subset\P(E)$ and a
direction $\C V\subset T_L\P(E)$ is identified to a direction (thus a
projective line $a$) through $L\in\P(E)$. $\varphi$ is, thus, a {\it
  homography} that associates to $A\in l$ (we identify $l$ with the
projective line $l/L\subset\P(E)$) the line $a\ni L$. As
$\varphi(L)=l$, we have, then, that three points
$(A,a),(B,b),(C,c)\in\beta^{(L,l)}$ belong to the same null-geodesic iff
\begin{equation}\label{birap}
(A,B:C,L)=(a,b:c,l), \end{equation}
i.e. the cross-ratio of the points $A,B,C,L\in l$ equals the
cross-ratio of the lines $a,b,c,l$ through $L$ (the dotted lines, and
their intersections with the lines $a,b,c$ are the points in the
integral $\alpha$-cone):
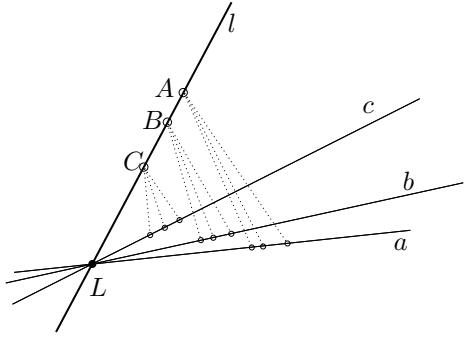
\begin{center}\begin{picture}(0,0)%
\epsfig{file=weyl/geodcp2.pstex}%
\end{picture}%
\setlength{\unitlength}{0.00035000in}%
\begingroup\makeatletter\ifx\SetFigFont\undefined
\def\x#1#2#3#4#5#6#7\relax{\def\x{#1#2#3#4#5#6}}%
\expandafter\x\fmtname xxxxxx\relax \def\y{splain}%
\ifx\x\y   
\gdef\SetFigFont#1#2#3{%
  \ifnum #1<17\tiny\else \ifnum #1<20\small\else
  \ifnum #1<24\normalsize\else \ifnum #1<29\large\else
  \ifnum #1<34\Large\else \ifnum #1<41\LARGE\else
     \huge\fi\fi\fi\fi\fi\fi
  \csname #3\endcsname}%
\else
\gdef\SetFigFont#1#2#3{\begingroup
  \count@#1\relax \ifnum 25<\count@\count@25\fi
  \def\x{\endgroup\@setsize\SetFigFont{#2pt}}%
  \expandafter\x
    \csname \romannumeral\the\count@ pt\expandafter\endcsname
    \csname @\romannumeral\the\count@ pt\endcsname
  \csname #3\endcsname}%
\fi
\fi\endgroup
\begin{picture}(6864,5001)(34,-4639)
\put(1322,-4038){\makebox(0,0)[lb]{\smash{\SetFigFont{10}{12.0}{rm}$L$}}}
\put(3392,-93){\makebox(0,0)[lb]{\smash{\SetFigFont{10}{12.0}{rm}$l$}}}
\put(1832,-2148){\makebox(0,0)[lb]{\smash{\SetFigFont{10}{12.0}{rm}$C$}}}
\put(5402,-1308){\makebox(0,0)[lb]{\smash{\SetFigFont{10}{12.0}{rm}$c$}}}
\put(6002,-2448){\makebox(0,0)[lb]{\smash{\SetFigFont{10}{12.0}{rm}$b$}}}
\put(5867,-3348){\makebox(0,0)[lb]{\smash{\SetFigFont{10}{12.0}{rm}$a$}}}
\put(2117,-1533){\makebox(0,0)[lb]{\smash{\SetFigFont{10}{12.0}{rm}$B$}}}
\put(2297,-1053){\makebox(0,0)[lb]{\smash{\SetFigFont{10}{12.0}{rm}$A$}}}
\end{picture}
\end{center}
We can now describe the null-geodesics passing through a point $(A,a)\in \mathbf{M}$ and
contained in a $\beta$-surface $\beta^(L,l)$, whose closure $\bar\beta$ is
isomorphic to $\cp1\times\cp1$: they coincide with the rational curves
in $\bar\beta$, containing $(A,a)$; except the ``horizontal''
($\bar\gamma^A$) and ``vertical'' ($\bar\gamma^a$) ones, all these
curves contain $(L,l)$:
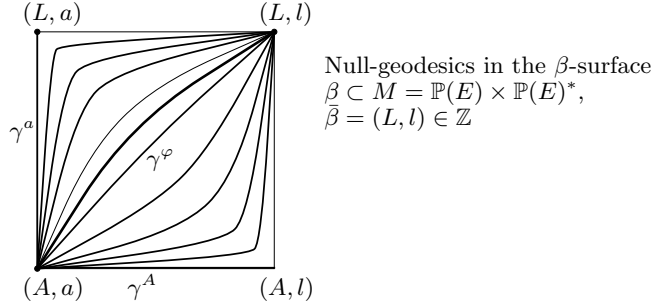
\begin{center}\input{weyl/geodbeta.pstex_t}\end{center}
We remark that, in the usual affine coordinates on
$\beta\simeq(\cp1\smallsetminus\{L\})\times
(\cp1\smallsetminus\{l\})\simeq\C^2$, these null-geodesics are the affine lines
containing $(A,a)$, thus the {\it projective structure} on $\beta$ is
(locally) isomorphic to a flat affine structure. We have seen, in
Section \ref{s.w6} (Corollary \ref{prflat}), that this is true for all
$\beta$-surfaces of a self-dual manifold. 
\smallskip

\subsection{The conformal structure of the complexification of $\cp2$}\label{ss.w86}
Let us study now the conformal structure of
$\mathbf{M}=\P(E)\times\P(E)^*\smallsetminus\mathcal{F}$ directly; actually {\bf M} has
a complex metric $g$. Let $(A,a)\in \mathbf{M}$, then $A$ is transverse to $a$,
thus we have the isomorphisms $E/a\simeq A$ and $E/A\simeq a$. Then,
a vector $(V,v)\in T_{(A,a)}\mathbf{M}$ is identified to a pair of
homomorphisms $V:A\rightarrow a$ and $v:a\rightarrow A$. Then the
metric $g$ is given by 
\begin{equation}
  \label{eq:metrcp2}
  g((V,v),(W,w)):=\mbox{tr}(v\circ W+w\circ V), \forall (V,v),(W,w)\in
  T_{(A,a)}\mathbf{M}.
\end{equation}
\obs {\bf (The real part). } Let $h$ be a hermitian metric on
$E$. Then we have a 
real-analytic embedding of $\mathbf{M}_0\simeq\P(E)$ into {\bf M}, given by:
$$\P(E)\ni A\longmapsto (A,A^\perp)\in
\P(E)\times\P(E)^*\smallsetminus\mathcal{F}.$$
A vector $(V,v)\in T_{(A,A^\perp)}\mathbf{M}$ is tangent to $\mathbf{M}_0$ iff
$$h(x,v(y))+h(V(x),y)=0, \;\forall x\in A,\;\forall y\in A^\perp.$$
Then one easily checks that $g((V,v),(W,w))=-2h(V,W), \,\forall
(V,v),(W,w)\in T_{(A,A^\perp)}\mathbf{M}_0$, hence, up to a constant, the
restriction of $g$ to $\mathbf{M}_0\simeq\cp2$ is the Fubini-Study metric of
$\cp2\simeq S^5/S^1$.
\medskip

An isotropic vector in {\bf M} is $(V,v)\in T_{(A,a)}\mathbf{M}$, such that $v\circ
V=0$, viewed as an endomorphism of $A$ (see above), or, equivalently,
such that 
\begin{equation}
  \label{eq:isotrcp2}
  \dim(A+V(A)\cap\ker v)>0.
\end{equation}
Let us see which is the limit of the isotropic cone in the points of
$\mathcal{F}$: from the relation above, it follows that the isotropic
cone in a point $x\in\mathcal{F}$ is 
$$C_x=\{(0,v)\in
T_x\mathcal{F}\}\cup \{(V,0)\in T_x\mathcal{F}\},$$ 
so the conformal
structure of {\bf M} is {\it singular} at the ``infinity'' $\mathcal{F}$.

\obs The situation $\mathcal{F}\subset\P(E)\times\P(E)^*$ is very
similar to the one treated in Section \ref{s.w7}, see also \cite{leb1}:
$\P(E)\times\P(E)^*$ has an Einstein self-dual metric $g$, singular at
the ``infinity'', and this Einstein structure yields a contact
structure on the twistor space $Z=\mathcal{F}$; the field of 2-planes
determined by this contact structure corresponds to the ``infinity''
$\mathcal{F}\subset \P(E)\times\P(E)^*$. But these planes do not admit
tangent rational curves, with normal bundle $\o 1\oplus\o 1$: the
conformal structure does not extend to the ``infinity'' (which is,
therefore, not a {\it conformal infinity}).
\medskip

\subsection{$\alpha$-planes and $\beta$-planes}\label{ss.w87}
We consider the isotropic planes in $T_{(A,a)}\mathbf{M}$ ($A\not\subset a$):
For a
fixed isotropic direction, represented by a generic vector $(V,v)\in
T_{(A,a)}\mathbf{M}$, the line $\ker v\subset a$ and the plane $V(A)+A\supset
A$ are fixed. The linear space of all vectors $(W,v)\in T_{(A,a)}\mathbf{M}$
satisfying
$$W(A)\subset A+V(A);\quad w|_{\ker v}=0$$
is isotropic and orthogonal to $(V,v)$: they form a $\beta$-plane. The $\alpha$-plane
$F^\alpha$ containing $(V,v)$ corresponds to the isotropic vectors
$(W,w)$, orthogonal to $(V,v)$, with $\ker w\neq \ker v$. 
As a plane transverse to all the $\beta$-planes (whose projection onto
$T_A\P(E)$ or $T_a\P(E)^*$ is never injective), $F^\alpha$ is
determined by a linear isomorphism $\varphi:T_A\P(E)\rightarrow
T_a\P(E)^*$, whose graph in $T_{(A,a)}\P(E)\times\P(E)^*$ is
$F^\alpha$; $\varphi$ induces the application
$\P\varphi:\P(a)\rightarrow\P(E/A)$ between the projective spaces of
$T_A\P(E)$, resp. $T_a\P(E)^*$. The plane $F^\alpha=F^\varphi$, the
graph of $\varphi$, is isotropic iff $V\subset\P\varphi (V)$, $\forall
V\in\P(a)$, i.e. $\P\varphi$ is the homography that sends a point $X$
in $a$ into the projective line through $A$ and $X$. We can extend
$\varphi$ to a projective isomorphism
$\varphi':\P(\C\oplus T_A\P(E))\rightarrow\P(\C\oplus T_a\P(E)^*)$: for example,
$\P(\C\oplus T_A\P(E))$ contains $T_A\P(E)$ as an affine open set. Then
$\varphi'$ is defined as follows:
\begin{eqnarray*}
  \varphi'|_{T_A\P(E)}&:=&\varphi\\
\varphi'|_{\P(T_A\P(E))}&:=&\P\varphi.
\end{eqnarray*}
Actually $\P(\C\oplus T_A\P(E))\simeq\P(E)$ and
$\P(\C\oplus T_a\P(E)^*)\simeq\P(E)^*$. We have then:
\begin{prop}
  A generic $\alpha$-plane $F^\alpha=F^\varphi$ in $T_{(A,a)}\mathbf{M}$ is
  the graph of a linear isomorphism $\varphi:T_A\P(E)\rightarrow
  T_a\P(E)^*$, which is determined by a projective isomorphism 
$$\varphi':\P(E)\rightarrow\P(E)^*,$$
such that $\varphi'(A)=a$ and $\varphi'(l)=l\cap a$, for all $l\supset A$.
\end{prop}
 
\subsection{Exponentials of $\alpha$-planes}\label{ss.w88} The exponential $\exp(F^\varphi)$
has an interpretation in terms of 
projective geometry: Each direction $\C(V,v)\subset F^\varphi$ is
determined by the point $\ker v$ in $a\subset\P(E)$ and the line
through $A$ and $\ker v$, and a homography $\phi^{(V,v)}$ from the
points $B$ of the projective line $A+\ker v$ to the space of lines $b$
through $\ker v$ (see next picture and the convention below). As this homography is the
restriction of $\varphi'$ to the appropriate spaces, it follows that
it is related to the homography $\phi^{(W,w)}$, where $\C(W,w)$ is
another direction in $F^\varphi$: the points $D:=b\cap c$,
$P:=a\cap (B+C)$ and $A$ are collinear:
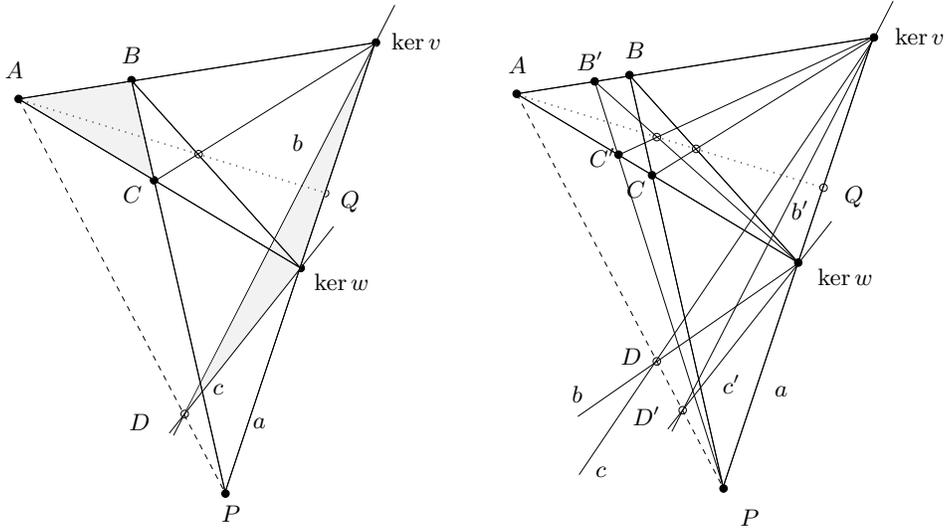
\begin{center}\input{weyl/alfacp2.pstex_t}\end{center}
Of course, this implies that $P$ determines a homography $\psi^P$
between the lines $A+\ker v$ and $A+\ker w$, such that $\psi^P(A)=A$
and $\psi^P(\ker v)=\ker w$. Then, for any other points $B'\in (A+\ker
v),\, C'=\psi^P(B)\in (A+\ker w)$, the lines $b'=\phi^{(V,v)}(B'),\,
c'=\phi^{(W,w)}(C')$ intersect on the line $(A+P)$ (see the right hand
side of the picture above).

{\noindent{\bf Convention. }} In the framework of plane projective
geometry, we identify a point in $\P(E)^*$ with a line in $\P(E)$ (we
note, for example $\ker v\in a$). The lines determined by the distinct
points $B$ and $C$ will be denoted by $(B+C)$ (thus $B,C\in(B+C)$).
\smallskip

The null-geodesic tangent to $(V,v)$ at $(A,a)$ is the set $\{(B,b)|B\in (A+\ker
v),\, b=\phi^{(V,v)}(B)\}$, and the null-geodesic tangent to $(W,w)$ is the
analogous set of the pairs $(C,c)$. Thus 
\begin{eqnarray*}
&\exp_{(A,a)}(F^\varphi)=\exp_{(A,a)}(F^\alpha)=\left\{(C,c)|C\in\P(E),\,
C\neq A,\right.&\\
&\left. c=((C+A)\cap a)+((A+P)\cap b^C)\right\}\cup\{(A,a)\},&
\end{eqnarray*}
where $B^C:=(A+\ker v)\cap (P+C)$, and $b^C:=\phi^{(V,v)}(B^C)$, as in
the picture above (where $B=B^C$, $b=b^C$). This gives the exponential
of the $\alpha$-plane determined by the isotropic vector $(V,v)$. We remark that
the point $(A,a)$ has a privileged position in
$\exp_{(A,a)}(F^\alpha)$: $(a\cap b)\in (A+B)$ $\forall (B,b)\in
\exp_{(A,a)}(F^\alpha)$; on the other hand $(b\cap c)\not\in (B+C)$ in
general (see the picture above), which means that the points $(B,b)$
and $(C,c)$ are not {\it null-separated} (i.e. they do not belong to
the same null-geodesic). That means that $\exp_{(A,a)}(F^\alpha)$ is not
totally isotropic, thus there is no $\alpha$-surface tangent to a generic $\alpha$-plane; not
surprising as the corresponding $\alpha$-cone is not flat (see
Section \ref{ss.w84}).

But there are $\alpha$-surfaces tangent to the two $\alpha$-planes $\{(V,0)|V\in T_A\P(E)\}$
and $\{(0,v)|v\in T_a\P(E)^*\}$: the ``slices'' $\{A\}\times\P(E)^*$
and $\P(E)\times\{a\}$. (It is easy to see that these planes are
isotropic, and that they are not $\beta$-planes, as these project on {\it lines}
in $T_A\P(E)$, resp. $T_a\P(E)^*$.)

Thus $\mathbf{M}=\P(E)\times\P(E)^*\smallsetminus\mathcal{F}$ is a conformal
self-dual manifold, not anti-self-dual, that admits $\alpha$-surfaces passing
through any point.

\bigskip
\begin{center}
{\sc {Centre de Math{\'e}matiques\\ UMR 7640 CNRS\\ Ecole Polytechnique
\\91128 Palaiseau cedex\\France}}\\e-mail: {\tt
belgun\@@math.polytechnique.fr}
\end{center}
\bigskip
\end{document}